\documentclass[11pt]{article}
\usepackage[utf8]{inputenc}

\usepackage{graphicx}
\usepackage{caption}
\usepackage{amsmath}
\usepackage{amssymb}
\usepackage{float}
\usepackage{hyperref}
\usepackage{amsthm}
\usepackage{xcolor}
\usepackage{comment}
\usepackage{amsfonts}
\usepackage{dsfont}
\usepackage{stmaryrd}
\usepackage{xcolor}
\usepackage{subcaption}

\usepackage{comment} \usepackage{siunitx}
\usepackage{mathtools} 
\mathtoolsset{showonlyrefs} 

\usepackage{lineno} 

\DeclareMathOperator{\R}{\mathbb{R}}

\def\Lc{{\cal L}}

\def\Gc{{\cal G}}
\def\Fc{{\cal F}}

\def\Nc{{\cal N}}
\def\Pc{{\cal P}}
\def\Sc{{\cal S}}
\def\Uc{{\cal U}}
\def\Xc{{\cal X}}

\def\Wc{{\cal W}}
\def\Zc{{\cal Z}}
\def \E{\mathbb{E}}

\def \P{\mathbb{P}}
\def \S{\mathbb{S}}

\def\1{{\bf 1}}

\def \N{\mathbb{N}}

\newcommand{\di}{\mathrm{d}}

\def \Sum{\displaystyle\sum}

\def \b1{\bf{1}}

\def\boX{{\boldsymbol X}}

\def\boW{{\boldsymbol W}}
\def\boZ{{\boldsymbol Z}}
\def\boZc{{\boldsymbol \Zc}}
\def\bog{{\boldsymbol \gamma}}
\def\boG{{\boldsymbol \Gamma}}

\def\bolx{{\boldsymbol x}}
\def\boly{{\boldsymbol y}}

\def\boz{{\boldsymbol z}}

\def \mfs{\mathfrak{s}}

\def \mrz{\mathrm{z}} 
 
\def \mfS{\mathfrak{S}} 
\def \bomrz{{\boldsymbol \mrz}}

\def\blue[#1]{{\color{blue}#1}}

\def\beqs{\begin{eqnarray*}}
\def\enqs{\end{eqnarray*}}
\def\beq{\begin{eqnarray}}
\def\enq{\end{eqnarray}}

\usepackage[algo2e,ruled,vlined]{algorithm2e} 
 \usepackage{algorithm}
 \usepackage{algorithmicx}
\usepackage{algpseudocode}
\usepackage{ marvosym }
\usepackage{csquotes}
\usepackage{hyperref}

\def \eps{\varepsilon}

\def\argmin{\mathop{\rm argmin}}
\def\argmin_#1{\underset{#1}{\mathrm{argmin\, }}}

\def \ep{\hbox{ }\hfill$\Box$}

\newtheorem{Theorem}{Theorem}[section]
\newtheorem{Definition}[Theorem]{Definition}

\newtheorem{Lemma}[Theorem]{Lemma}

\newtheorem{Remark}[Theorem]{Remark}
\newtheorem{Example}[Theorem]{Example}

\numberwithin{equation}{section} 

\usepackage[backend=bibtex, style=alphabetic, maxbibnames=99]{biblatex}

\def \trans{^{\scriptscriptstyle{\intercal}}}

\begin{document}
\title{DeepSets and their derivative networks \\ 
for solving symmetric PDEs
\thanks{This work  was supported by FiME (Finance for Energy Market Research Centre) and
the ``Finance et D\'eveloppement Durable - Approches Quantitatives'' EDF - CACIB Chair.}
}

\author{
Maximilien \textsc{Germain}\footnote{EDF R\&D, and  LPSM, Université de Paris  \sf \href{mailto:Maximilien.Germain at edf.fr}{mgermain at lpsm.paris}} 
\and Mathieu \textsc{Lauri\`ere} \footnote{ORFE, Princeton University  \sf \href{mailto:lauriere at princeton.edu}{lauriere at princeton.edu}} 
\and  Huyên \textsc{Pham} \footnote{LPSM, Université de Paris, and  FiME,  and CREST ENSAE \sf \href{mailto:pham at lpsm.paris}{pham at lpsm.paris}} 
\and  Xavier \textsc{Warin} \footnote{EDF R\&D, and FiME \sf \href{mailto:xavier.warin at  edf.fr}{xavier.warin at edf.fr}} 
}

\date{\today\\ {\it to appear in Journal of Scientific Computing}}

\maketitle              
\begin{abstract}
Machine learning methods for solving nonlinear partial differential equations (PDEs) are hot topical issues, and different algorithms proposed in the literature  show efficient numerical approximation  in high dimension.
In this paper, we introduce a class of  PDEs that are invariant to permutations, and called {\it symmetric} PDEs. Such pro\-blems are widespread, ranging from cosmology to quantum mechanics, and option pricing/hedging in multi-asset market with exchangeable payoff.  Our main application comes actually from the particles approximation of mean-field control problems.  We design deep lear\-ning algorithms based on certain types of neural networks, named PointNet and DeepSet (and their associated derivative networks), for computing simultaneously an approximation of the solution and its gradient to symmetric PDEs. 
We illustrate the performance and accuracy  of the PointNet/DeepSet networks compared to classical feedforward ones, and provide several numerical results of our algorithm for the examples of a mean-field systemic risk,  mean-variance problem and a min/max linear quadratic 
McKean-Vlasov control problem. 
\end{abstract}

\vspace{5mm}


\section{Introduction}

The numerical resolution of partial differential equations (PDEs) in high dimension  is a major challenge in various areas of science, engineering, and finance. PDEs that appear in the applications are often non linear and of very high dimension (number of particles in physics, 
number of agents in large population control problems,  number of assets and factors in financial markets, etc), and are subject to the so-called curse of dimensionality, which makes infeasible the implementation of  classical grid methods and Monte-Carlo approaches. 

\vspace{1mm}

A breakthrough with deep learning based-algorithms has been made in the last five years towards this computational challenge, and we mention the recent survey papers by \cite{bechutjenkuc20} and \cite{GPW21}.  
The main interest in the use of machine learning techniques for PDEs is the ability of deep neural networks to efficiently represent high dimensional functions without using spatial grids, and with no curse of dimensionality  (see e.g. \cite{hutetal20}).  
Although the use of neural networks for solving PDEs is not new,  the approach has been successfully revived with new ideas and directions.  Moreover, recently developed open source libraries like Tensorflow 
and  Pytorch  offer an accessible framework to implement these algorithms. 

\vspace{1mm}

In this paper, we introduce a class  of PDEs  that are invariant by permutation, and called here  {\it symmetric} PDEs.  Such  PDEs occur naturally in the modelling of systems dea\-ling with sets that are invariant by permutation of their elements. Applications range from  models in general relativity and cosmology, to  quantum mechanics and  chemistry, see  e.g. 
\cite{vanwyk08}, \cite{smu11}.  Symmetric PDEs  also appear in the pricing/hedging of basket option and  options on the maximum of multiple assets.  Our main motivation for introducing this general class of symmetric PDEs comes from the  control of 
large population of interacting indistinguishable agents, which leads in the asymptotic regime of infinite population to the  theory of mean-field games (MFG) and mean-field type control, also called McKean-Vlasov (MKV)  control. These topics  have 
attracted an increasing and large interest since the seminal papers  \cite{laslio07} and \cite{huacaimal06} with important mathematical developments and numerous applications in various fields over the last decade. 
We refer to the two-volume monographs \cite{cardel19}-\cite{cardel2} for an exhaustive exposition of this research domain, where it is known that the solution to MFG or  MKV control problem are characterized  in terms of a Master equation or a Bellman equation, which are PDEs in the Wasserstein space of probability measures.  It  turns out that the finite-dimensional approximation of these equations  are formulated as  {\it symmetric} non linear PDEs, and the convergence of this approximation has been recently obtained in  \cite{ganmayswi20}, and  \cite{GPW21b} 
(for a rate of convergence), see also  \cite{lac17} and  \cite{dje20}.  

\vspace{1mm}

Symmetric PDEs are often in very high dimension, typically of the order  of one thousand in the case of particles approximation of Master and Bellman equations, and it is tempting to apply machine learning algorithms in this framework. For that purpose, we shall furthermore exploit the symmetric structure of the PDEs in order to design deep learning-based algorithms with a suitable class of neural networks. We first observe that the solution to symmetric PDEs is invariant by permutation (also called {\it exchangeable}),   
and we shall then consider a class of symmetric neural networks, named PointNet and DeepSets, aiming to approximate our solution. 
Such class of neural networks has been recently introduced in the machine learning community, see \cite{pointnet}, \cite{deepsets}, \cite{bloteh20},  for dealing with tasks invol\-ving some invariant data sets,  and it turns out that they provide much better  accuracy than classical feedforward neural networks (NN in short) in the approximation of symmetric functions. Indeed, feedforward NN have too many degrees of freedom, and the optimization over parameters in (stochastic) gradient descent algorithm may be trapped away in the approximation of a symmetric function, as illustrated in several examples and comparison tests 
presented in this paper. 
We shall also introduce different  classes of derivative symmetric network, named DeepDerSet and AD-DeepSet, for the approximation of the gradient of the solution to symmetric PDEs.

\vspace{1mm}

By relying on the class of symmetric NN, and their derivative networks, we next adapt the deep backward dynamic programming scheme \cite{hure2020deep}, \cite{phawarger20} for numerically sol\-ving  symmetric PDEs, i.e., finding approximations of the solution and its gradient. We also explain 
in the case of mean-field control problem how our scheme provides an approximation for the solution to a Bellman equation in the Wasserstein space of probability measures. This yields alternative deep learning schemes for mean-field control problems  to the ones recently designed in \cite{germicwar19}, \cite{carlau19}, \cite{fouquezhang19}, or \cite{Ruthotto9183}. 
We test our algorithms on several examples arising from different McKean-Vlasov control problem,  for which we have explicit or benchmarked solutions: a systemic risk model as in \cite{carfousun}, the classical mean-variance, i.e., Markowitz portfolio allocation problem, and  a 
min/max linear quadratic mean-field control problem as in \cite{salhab2015a}. 

\vspace{1mm} 
  
 \noindent {\bf Outline of the paper.}  The rest of the paper is organized as follows. We introduce in Section \ref{secsymPDE} the class of {\it symmetric} PDEs with some examples, and  show  exchangeability properties of the solution and its gradient to such PDEs. 
 Section \ref{secsymNN} is devoted to the exposition of the class of symmetric neural networks, as well as its derivative networks, and we provide several comparison tests with respect to classical feedforward NN.  We describe in Section \ref{secalgo}  
 the deep learning schemes for solving symmetric PDEs, and finally provide several numerical examples in Section \ref{secnum}. 
 
 \vspace{1mm}

\noindent {\bf Notations.}   Given $N$ $\in$ $\N^*$, $\Xc^N$ denotes the set of all elements $\bolx$ $=$ $(x_i)_{i\in\llbracket 1,N\rrbracket}$ with coefficients $x_i$ valued in $\Xc$ and $\llbracket 1,N \rrbracket=\{1,\cdots,N\}$. When $\Xc$ $=$ $\R^d$, one usually identifies  $(\R^d)^N$ with $\R^{d\times N}$  the set of $d\times N$-matrices with real-valued coefficients.
$\S^N(\Xc)$ is  the set of $N\times N$-symmetric matrices with coefficients valued in $\Xc$, 
and is simply denoted by  $\S^N$ when $\Xc$ $=$ $\R$.   For a real-valued $C^2$ function $\varphi$ defined on $(\R^d)^N$, its gradient $D\varphi(\bolx)$ $=$ $(D_{x_i}\varphi(\bolx))_{i\in\llbracket 1,N\rrbracket}$ is valued in $(\R^d)^N$, while its Hessian 
$D^2\varphi(\bolx)$ $=$ $(D_{x_ix_j}^2\varphi(\bolx))_{i,j\in\llbracket 1,N \rrbracket}$ is  valued in $\S^N(\S^d)$.

We denote by $\mfS_N$ the set of permutations  on $\{1,\ldots,N\}$. For any $\bolx$ $=$  $(x_i)_{i\in\llbracket 1,\N \rrbracket}$   $\in$ $\Xc^N$, $\pi$ $\in$ $\mfS_N$, we denote by   
$\pi[\bolx]$ $=$ $(x_{\pi(i)})_{i\in\llbracket 1,\N \rrbracket}$  $\in$ $\Xc^N$.  For any $\boG$ $=$ $(\Gamma_{ij})_{i,j\in\llbracket 1,\N \rrbracket}$ $\in$ $\S^{N}(\Xc)$,  we denote by $\pi[\boG]$ $=$ 
$(\Gamma_{\pi(i)\pi(j)})_{i,j\in\llbracket 1,N \rrbracket}$ $\in$ $\S^{N}(\Xc)$.

We say that a function $\varphi$ defined on $\Xc^N$ is exchangeable to the order $N$  on $\Xc$  if it is invariant by permutation, i.e., 
$\varphi(\bolx)$ $=$ $\varphi(\pi[\bolx])$, for any $\bolx$ $\in$ $\Xc^N$, and $\pi$ $\in$ $\mfS_N$. We may  simply say exchangeable when it is clear from the context. 
The notations $0_d, 1_d$ refer respectively to $d$-dimensional vectors full of $0$ and $1$. With two vectors $a,b\in\R^d$, the notation $a.b=\sum_{i=1}^d a_i b_i$ refers to the canonical scalar product.
 
\section{Symmetric PDEs} \label{secsymPDE} 

We consider a so-called {\it symmetric} class of parabolic second-order partial differential equations (PDEs): 
\begin{equation} \label{symPDE}
\left\{
\begin{array}{rcl}
\partial_t v  +  F(t,\bolx,v,D_{\bolx} v,D_\bolx^2v) &=& 0, \quad \quad  (t,\bolx) \in [0,T) \times (\R^d)^N  \\
v(T,\bolx) &=&  G(\bolx), \quad \bolx \in (\R^d)^N, 
\end{array}
\right.
\end{equation}
where $F$  is a  real-valued function defined on $[0,T]\times(\R^d)^N\times\R\times(\R^d)^N\times\S^{N}(\S^d)$, $G$ is defined on  $(\R^d)^N$,   and satisfying the permutation-invariance condition:

\vspace{2mm}

\noindent {\bf (HI)} \quad For  any  $t$ $\in$ $[0,T]$, $\bolx$ $\in$ $(\R^d)^N$, $y$ $\in$ $\R$, $\boz$ $\in$ $(\R^d)^N$, $\bog$ $\in$ $\S^{N}(\S^d)$, 
\beqs
F(t,\bolx,y,\boz,\bog) &=& F(t,\pi[\bolx],y,\pi[\boz],\pi[\bog]) \\
G(\bolx) &=& G(\pi[\bolx]), \quad \forall \pi \in \mfS_N.   
\enqs
We assume that PDE \eqref{symPDE} is well-posed in the sense that there exists a unique classical solution satisfying a suitable growth condition. 

\vspace{1mm}

We list below some examples of symmetric PDEs in the form \eqref{symPDE}.  We start with an example of pricing in a ``symmetric" multi-asset model.

 \begin{Example}[{\it Multi-asset pricing}] \label{exmultiasset} 
 {\rm
 Let us consider a model with $N$  risky assets of price process $\boX$ $=$ $(X^1,\ldots,X^N)$ governed by 
 \beqs
 dX_t^i &=& \tilde b_i(\boX_t) dt +  \sum_{j=1}^N \sigma_{ij}(\boX_t) dW_t^j,
 \enqs
 where the diffusion coefficients satisfy the property: for all $\pi$ $\in$ $\mfS_N$,
 \begin{equation} \label{symsig} 
 \sigma_{ij}(\pi[\bolx]) \; = \;  \sigma_{\pi(i)\pi(j)}(\bolx), \quad \bolx  = (x_i)_{i\in\llbracket 1,N\rrbracket}, \; i,j=1,\ldots,N. 
 \end{equation}
 Notice that $\tilde b_i$ is the drift of the asset price under the historical probability measure, and does not appear in the pricing equations below. 
 The symmetry condition \eqref{symsig}  is satisfied for example when $\sigma_{ii}(\bolx)$ $=$ $\sigma(\bolx)$, and $\sigma_{ij}(\bolx)$ $=$ $\tilde\sigma(\bolx)$, $i,j=1,\ldots,N$, $i\neq j$, with $\sigma, \tilde\sigma$ exchangeable functions. Another example 
 is when $\sigma_{ii}(\bolx)$ $=$ $\sigma(x_i)$, and $\sigma_{ij}(\bolx)$ $=$ $\tilde\vartheta(x_i)\bar\vartheta(x_j)$, $i,j=1,\ldots,N$, $i\neq j$, for some functions $\sigma$, $\tilde\vartheta$, $\bar\vartheta$ defined on $\R$, 
 which means that all the assets have the same marginal volatility coefficient, and the correlation function between any pair of assets is identical. We 
 consider an option of maturity $T$ with payoff $G(X_T^1,\ldots,X_T^N)$, where $G$ is an exchangeable function, for example: 
 \beqs
 G(\bolx) &=& 
 \left\{
 \begin{array}{ll} 
 \big( \max(x_1,\ldots,x_N) -K\big)_+,  & \mbox{ (call on max) }  \\
\big( \sum_{i=1}^N x_i - K \big)_+, & \mbox{ (call on sum)},  \\
\sum_{i=1}^N 1_{x_i \geq K}, &  \mbox{ (sum of binary options)},
 \end{array}
 \right.
 \enqs
 for $\bolx$ $=$ $(x_1,\ldots,x_N) \in \R^N$.  
 In a frictionless market with constant interest rate $r$, the option price $(t,\bolx)$ $\in$ $[0,T]\times\R^N$  $\mapsto$ $v(t,\bolx)$ satisfies a linear PDE \eqref{symPDE} with terminal condition given by the exchangeable function $G$ and 
 \beqs
 F(t,\bolx,y,\boz,\bog) &=& -ry + r \sum_{i=1}^N  x_i z_i  + \frac{1}{2} \sum_{i,j=1}^N \sigma_{ij}^2(\bolx) \gamma_{ij}, 
 \enqs
 for  $t$ $\in$ $[0,T]$, $\bolx$ $=$ $(x_i)_{i\in\llbracket 1,N \rrbracket}$  $\in$ $\R^N$, $y$ $\in$ $\R$, $\boz$ $=$ $(z_i)_{i\in\llbracket 1,\N \rrbracket}$ $\in$ $\R^N$, and $\bog$ $=$ $(\gamma_{ij})_{i,j\in\llbracket 1,N \rrbracket}$ $\in$ $\S^{N}$.  
 In the case of counterparty risk, the pricing of CVA leads to a quasi-linear PDE \eqref{symPDE} with $F$ in the form (see \cite{PHL12} for the details of the PDE derivation): 
 \beqs
 F(t,\bolx,y,\boz,\bog) &=&  \beta(y^+ - y) + r \sum_{i=1}^N x_i  z_i    + \frac{1}{2} \sum_{i,j=1}^N \sigma_{ij}^2(\bolx) \gamma_{ij}, 
 \enqs
where $\beta$ $>$ $0$ is the intensity of default.  Another case of non-linearity occurs when lending rate $r>0$ is smaller than borrowing rate $R>0$, which leads to a super-replication price solution to  \eqref{symPDE} with $F$ given by 
 \beqs
 F(t,\bolx,y,\boz,\bog) &=&  \sup_{b \in [r,R]} \big[ -by + b \sum_{i=1}^N  x_i z_i \big]  + \frac{1}{2} \sum_{i,j=1}^N \sigma_{ij}^2(\bolx) \gamma_{ij}.
 \enqs
In the above three cases, and under \eqref{symsig}, the generator function $F$ clearly satisfies the permutation-invariance condition in  {\bf (HI)}. 
 }
 \ep
 \end{Example}

\vspace{1mm}

The second example is actually our main motivation for considering symmetric  PDEs, and comes from mean-field models. 

\begin{Example}[{\it McKean-Vlasov control problem with common noise}] \label{exMKV} 
{\rm Let us consider $N$ interacting indistinguishable agents with controlled  state process $\boX$ $=$ $(X^1,\ldots,X^N)$  valued in $(\R^d)^N$, and  driven by  
\begin{align}
dX_t^i & = \;  \beta(t,X_t^i,\bar\mu(\boX_t),\alpha_t^i) dt + \sigma(t,X_t^i,\bar\mu(\boX_t),\alpha_t^i) dW_t^i  \nonumber \\
&  \quad \quad +  \; \sigma_0(t,X_t^i,\bar\mu(\boX_t)) dW_t^0,  \quad \quad 0 \leq t \leq T, \; i=1,\ldots,N, \nonumber 
\end{align}
where $\bolx$ $=$ $(x_i)_{i\in\llbracket 1,N\rrbracket}$ $\mapsto$ $\bar\mu(\bolx)$ $=$ $\frac{1}{N}\sum_{i=1}^N \delta_{x_i}$ is the empirical measure (exchangeable) function, 
$W^i$, $i$ $=$ $1,\ldots,N$, are independent Brownian motions representing idiosyncratic noises, and $W^0$ is a Brownian motion independent of 
$\boW$ $=$ $(W^i)_{i\in\llbracket 1,N\rrbracket}$, representing a common noise.  Moreover, $\alpha^i$ is a control process (valued in some Polish space $A$) applied by the agent $i$ who follows in a cooperative equilibrium a social planner aiming to minimize a social cost in the form 
\beqs
J(\alpha^1,\ldots,\alpha^N) &=& \frac{1}{N} \sum_{i=1}^N  \E \Big[ \int_0^T e^{-rt} f(X_t^i,\bar\mu(\boX_t),\alpha_t^i) dt + e^{-rT} g(X_T^i,\bar\mu(\boX_T)) \Big]. 
\enqs
The Bellman equation to this $N$-cooperative agents control problem is in the form \eqref{symPDE} with a Hamiltonian function $F$ given by 
\beqs
F(t,\bolx,y,\boz,\bog) &=& \sum_{i=1}^N \inf_{a \in A} \big[  \beta(t,x_i,\bar\mu(\bolx),a).z_i + \frac{1}{2} {\rm tr}\big(\Sigma(t,x_i,\bar\mu(\bolx),a)\gamma_{ii} \big) + \frac{1}{N}  f(x_i,\bar\mu(\bolx),a) \big] \\
& & \quad \quad \quad \quad + \; \frac{1}{2} \sum_{i \neq j } {\rm tr}\big(\sigma_0(t,x_i,\bar\mu(\bolx))\sigma_0\trans(t,x_j,\bar\mu(\bolx))\gamma_{ij} \big) \; - ry, 
\enqs
where $\Sigma$ $=$ $\sigma\sigma\trans+\sigma_0\sigma_0\trans$, and a terminal condition given by 
\beqs
G(\bolx) &=& \frac{1}{N} \sum_{i=1}^N g(x_i,\bar\mu(\bolx)). 
\enqs
Such functions $F$ and $G$ clearly satisfy condition {\bf (HI)}.  Let us point out that in the limiting regime when the number $N$ of agents goes to infinity, it is proved in \cite{lac17}, \cite{dje20}, \cite{ganmayswi20}   that the solution to this cooperative-agents problem converges to the McKean-Vlasov 
control problem with state process $X$ $=$ $X^\alpha$ of dynamics
\begin{align} \label{diffXMKV} 
dX_t &= \beta(t,X_t,\P^0_{X_t},\alpha_t) dt + \sigma(t,X_t,\P^0_{X_t},\alpha_t) dW_t +  \sigma_0(t,X_t,\P^0_{X_t}) dW_t^0, 
\end{align}
and cost functional 
\beqs
J_{_{MKV}}(\alpha) &=& \E\Big[ \int_0^T e^{-rt} f(X_t,\P^0_{X_t},\alpha_t) dt + e^{-rT} g(X_T,\P^0_{X_T}) \Big]. 
\enqs
(Here $\P^0_{X_t}$ denotes the conditional law of $X_t$ given the common noise $W^0$). 
Moreover, the corresponding Bellman equation in the Wasserstein space of square-integrable probability measures $\Pc_{2}(\R^d)$ is given by (see \cite{phawei17}) 
\begin{equation} \label{meanPDE}
\left\{
\begin{array}{rcl}
\partial_t v  +  \Fc(t,\mu,v,\partial_\mu v, \partial_x\partial_\mu v,\partial_\mu^2 v) &=& 0, \quad \quad  (t,\mu) \in [0,T) \times \Pc_2(\R^d)  \\
v(T,\mu) &=&  \Gc(\mu), \quad \mu \in \Pc_2(\R^d), 
\end{array}
\right.
\end{equation}
where $\partial_\mu \varphi(\mu)(.)$ $:$ $\R^d$ $\rightarrow$ $\R^d$, $\partial_x\partial_\mu \varphi(\mu)(.)$  $:$ $\R^d$ $\rightarrow$ $\S^{d}$, $\partial_\mu^2 \varphi(\mu)(.,.)$ $:$ $\R^d\times\R^d$ $\rightarrow$ $\S^d$, are the $L$-derivatives of a function $\varphi$ on $\Pc_{2}(\R^d)$ (see \cite{cardel19}) 
and 
\beqs
\Fc(t,\mu,y,Z(.),\Gamma(.),\Gamma_0(.,.)) &=& -ry + \int_{\R^d} h(t,x,\mu,Z(x),\Gamma(x)) \mu(dx) \\ 
& & \;  + \;  \int_{\R^d\times\R^d}  \frac{1}{2} {\rm tr}\big( \sigma_0(t,x,\mu)\sigma_0\trans(t,x',\mu) \Gamma_0(x,x') \big) \mu(dx)\mu(dx'),  \\
\Gc(\mu) &=& \int_{\R^d} g(x,\mu) \mu(dx),  
\enqs
with 
\beqs
h(t,x,\mu,z,\gamma) &=&  \inf_{a\in A} \Big[ \beta(t,x,\mu,a).z +  \frac{1}{2} {\rm tr}\big(\Sigma(t,x,\mu,a)\gamma \big)  + \;  f(x,\mu,a) \Big]. 
\enqs
}
\ep
\end{Example}

\vspace{1mm}

We end this section by showing some exchangeability properties of the solution to the symmetric PDE \eqref{symPDE}.   
Let us introduce the notion of $D$-exchangeability where $D$ stands for derivative.

 \begin{Definition}
 A function $(\bolx,x)$ $\in$ $(\R^d)^N\times\R^d$ $\mapsto$ $\mrz(\bolx,x)$ $\in$ $\Xc$  is $D$-exchangeable if for any fixed $x$ $\in$ $\R^d$,  $\mrz(.,x)$ is exchangeable. Given a $D$-exchangeable function $\mrz$, we denote by 
$\bomrz$  the function defined on $(\R^d)^N$ by  $\bomrz(\bolx)$ $=$ $(\mrz(\bolx,x_i))_{i\in\llbracket 1,N\rrbracket}$  $\in$  $\Xc^N$. 
 \end{Definition}

This definition is actually motivated by the exchangeability property of the solution to the PDE  \eqref{symPDE}, and by a structural property on the gradient of an  exchangeable function that is differentiable. 

\begin{Lemma} \label{lemvsym}
The solution $v$ to the PDE \eqref{symPDE} with $F$ and G satisfying {\bf (HI)} is exchangeable, i.e.,   for all $\pi$ $\in$ $\mfS_N$, 
\beqs
v(t,\bolx) &=& v(t,\pi[\bolx]),  \quad  (t,\bolx) \in [0,T]\times(\R^d)^N.  
\enqs
\end{Lemma}
\noindent {\bf Proof.} Let $\pi$ $\in$ $\mfS_N$. We set $v_\pi(t,\bolx)$ $=$ $v(t,\pi[\bolx])$, and observe that  $\partial_t v(t,\pi[\bolx])$ $=$ $\partial_t v_\pi(t,\bolx)$, while 
\beqs
D_{\bolx} v(t,\pi[\bolx]) \; = \; \pi[ D_\bolx v_\pi(t,\bolx) ], \quad D^2_{\bolx} v(t,\pi[\bolx]) \; = \; \pi[ D_\bolx^2 v_\pi(t,\bolx) ]. 
\enqs
By writing the PDE \eqref{symPDE} at $(t,\pi[\bolx])$, it follows under {\bf(HI)} that $v_\pi$ satisfies 
\begin{equation} 
\left\{
\begin{array}{rcl}
\partial_t v_\pi  +  F(t,\bolx,v_\pi,D_{\bolx} v_\pi,D_\bolx^2v_\pi) &=& 0, \quad \quad  (t,\bolx) \in [0,T) \in (\R^d)^N  \\
v_\pi(T,\bolx) &=&  G(\bolx), \quad \bolx \in (\R^d)^N. 
\end{array}
\right.
\end{equation}
By uniqueness of the solution to PDE \eqref{symPDE}, we conclude that $v_\pi$ $=$ $v$, i.e., the exchan\-geability property of $v$. 
\ep

\vspace{1mm}

\begin{Lemma} \label{lemvNderiv} 
Let $w$ be an exchangeable, and differentiable function on $(\R^d)^N$.  Then there exists a $D$-exchangeable function $\mrz$ such that 
\begin{align} \label{excz1}
D_{x_i} w(\bolx) &= \; \mrz(\bolx,x_i), \quad i=1,\ldots,N, 
\end{align}
for all $\bolx$ $=$ $(x_i)_{\in\llbracket 1,N\rrbracket}$ $\in$ $(\R^d)^N$, i.e.,  $Dw$ $=$ $\bomrz$. 
\end{Lemma}  
\noindent  {\bf Proof.} Since $w$ is exchangeable, it is clear that for fixed $i$ $\in$ $\llbracket 1,N\rrbracket$, and $x_i$ $\in$ $\R^d$,
\begin{align}
\bolx_{-i} := (x_j)_{j\neq i} \in (\R^d)^{N-1} & \mapsto \;  D_{x_i} w(x_1,\ldots,x_{i-1},x_i,x_{i+1},\ldots,x_N) \quad \mbox{ is exchangeable}, 
\end{align} 
and  we shall then write: 
\begin{align}
\mrz^i(\bolx_{-i},x) & := \; D_{x_i} w (x_1,\ldots,x_{i-1},x,x_{i+1},\ldots,x_N), \quad x \in \R^d.  
\end{align}
By exchangeability of $w$, we also note that 
\begin{align} \label{exil1} 
\mrz^i(\bolx_{-i},x) &= \; \mrz^\ell(\bolx_{-i},x), \quad \forall  i,\ell  \in \llbracket1, N \rrbracket.  
\end{align}
Let us now define the function $\mrz$ on $(\R^d)^N\times\R^d$ by:
\begin{align} \label{defzk} 
\mrz(\bolx,x) &:= \; \frac{1}{N}  \sum_{p=0}^{N-1} (-1)^p  \sum_{1\leq i_1<\ldots<i_{p+1}\leq N} 
\sum_{\ell=1}^N \mrz^\ell((\underbrace{x,\ldots,x}_{p\; \mbox{\scriptsize{times}}},x_j)_{j\neq i_1,\ldots,i_{p+1}},x).  
\end{align}
for $\bolx$ $=$ $(x_1,\ldots,x_N)$ $\in$ $(\R^d)^N$, and $x$ $\in$ $\R^d$,  $(\underbrace{x,\ldots,x}_{p\; \mbox{\scriptsize{times}}},(x_j)_{j\neq i_1,\ldots,i_{p+1}})$ is the vector in $(\R^d)^{N-1}$ consisting of $p$ components $x$, and the $N-p-1$ components 
$x_j$, for $j$ $\neq$ $i_1,\ldots,i_{p+1}$.
By construction, it is clear  that for fixed $x$ $\in$ $\R^d$, $\mrz(.,x)$ is exchangeable, i.e., $\mrz$ is a $D$-exchangeable function.   
Let us now show \eqref{excz1}, i.e.,  that for  fixed  $\bolx$ $=$ $(x_i)_{i\in\llbracket 1,N\rrbracket}$ $\in$ $(\R^d)^N$,  
\begin{align}
\mrz(\bolx,x_i) &= \mrz^i(\bolx_{-i},x_i), \quad i \in   \llbracket 1, N\rrbracket. 
\end{align} 
It suffices to check this property for $i$ $=$ $1$.  We set for $p$ $=$ $0,\ldots,N-1$:
\begin{align}
S^p & := \;   \sum_{1\leq i_1<\ldots<i_{p+1}\leq N} \sum_{\ell=1}^N \mrz^\ell(x_1,(\underbrace{x_1,\ldots,x_1}_{p\; \mbox{\scriptsize{times}}},x_j)_{j\neq i_1,\ldots,i_{p+1}}).
\end{align}
and see that 
\begin{equation}
\begin{cases}
S^0 &= \; \Sum_{\ell=1}^N \mrz^\ell(\bolx_{-1},x_1)  + \Sum_{i_1=2}^N \sum_{\ell=1}^N \mrz^\ell((x_1,x_j)_{j\neq 1,i_1},x_1) \\
S^1 & = \; \Sum_{i_2=2}^N \Sum_{\ell=1}^N \mrz^\ell((x_1,x_j)_{j\neq 1,i_2},x_1)  + \Sum_{2\leq i_1 < i_2\leq N}  \Sum_{\ell=1}^N \mrz^\ell((x_1,x_j)_{j\neq i_1,i_2},x_1) \\
& \vdots \\
S^{N-2} &= \;  \Sum_{2\leq i_2 <\ldots< i_{N-1} \leq N} \Sum_{\ell=1}^N \mrz^\ell((\underbrace{x_1,\ldots,x_1}_{N-2\; \mbox{\scriptsize{times}}},x_j)_{j\neq 1,i_2,\ldots,i_{N-1}},x_1) 
+  \Sum_{\ell=1}^N \mrz^\ell((x_1,\ldots,x_1),x_1) \\
S^{N-1} &= \; \Sum_{\ell=1}^N \mrz^\ell((x_1,\ldots,x_1),x_1). 
\end{cases}
\end{equation}
The telescopic sum then yields 
\begin{align}
\mrz((x_1,\ldots,x_N),x_1) &= \; \frac{1}{N}  \sum_{p=0}^{N-1} (-1)^p  S^p \; = \;  \frac{1}{N} \sum_{\ell=1}^N \mrz^\ell((x_j)_{j\neq 1},x_1) \; = \; \mrz^1((x_j)_{j\neq 1},x_1), 
\end{align}
where the last equality follows from \eqref{exil1}.  This shows the property \eqref{excz1}. 
 \ep

\section{Symmetric neural networks} \label{secsymNN}

\subsection{DeepSets and variants}

In view of Lemma \ref{lemvsym} and \ref{lemvNderiv}, we shall consider a class of neural networks (NN in short) that satisfy the exchangeability and $D$-exchangeability properties for approximating the solution (and its gradient) to the PDE \eqref{symPDE}.

We denote by 
\begin{align}
\Lc_{d_1,d_2}^{\rho} &=\;  \Big\{ \phi : \R^{d_1} \rightarrow \R^{d_2}: \exists \;  (\Wc,\beta) \in   \R^{d_2\times d_1} \times \R^{d_2}, \; 
\phi(x) \; = \: \rho( \Wc x + \beta) \; \Big\}, 
\end{align}
the set of layer  functions with input dimension $d_1$, output dimension $d_2$, and activation function $\rho$ $:$ $\R$ $\rightarrow$ $\R$.  Here, the activation is applied component-wise, i.e., 
$\rho(x_1,\ldots,x_{d_2})$ $=$ $\big(\rho(x_1),\ldots,\rho(x_{d_2})\big)$, to the result of the affine map $x$ $\in$ $\R^{d_1}$ $\mapsto$ $\Wc x + \beta$ $\in$ $\R^{d_2}$, 
with a matrix $\Wc$ called weight, and vector $\beta$ called bias.    
Standard examples of activation functions are  the sigmoid, the ReLU, the Elu (see \cite{elu}), or $\tanh$. 
When $\rho$ is the identity function, we simply write $\Lc_{d_1,d_2}$.   

\vspace{1mm}

We then define 
\begin{align}
\Nc^\rho_{d_0,\ell,m,k} &= \;  \Big\{ \varphi : \R^{d_0}  \rightarrow \R^k: \exists   \phi_0 \in  \Lc^\rho_{d_0,m}, \; \exists \phi_i \in \Lc^\rho_{m,m}, i=1,\ldots,\ell-1, \; \exists  \phi_\ell \in \Lc_{m,k},  \\
& \hspace{4cm}   \varphi  \; = \;  \phi_\ell \circ \phi_{\ell -1} \circ \cdots \circ \phi_0 \Big\},
\end{align}
as the set of feedforward (or artificial) neural networks with input layer dimension $d_0$, output layer dimension $k$, and $\ell$  hidden layers with $m$ neurons (or units).  
These numbers $d_0,\ell,m$, and the activation function $\rho$,  form the architecture of the network. 
When $\ell$ $=$ $1$, one usually refers to shallow neural networks,  as opposed to deep neural networks which have several hidden layers.

\vspace{1mm}

\begin{figure}[H]
    \centering
    \includegraphics[width=10cm]{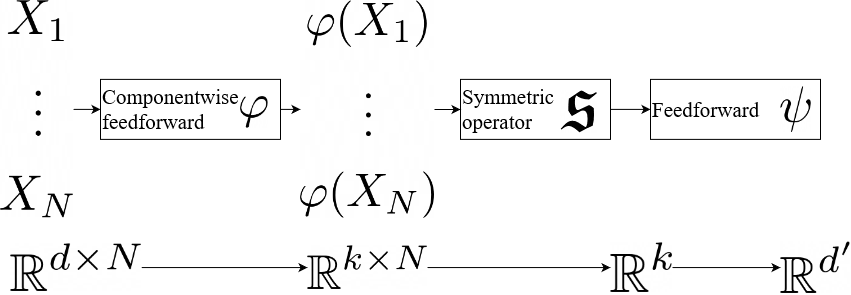}
     \caption{Architecture of a symmetric neural network.
     }
\label{fig: symNN}
\end{figure}

A {\it symmetric neural network} function, denoted  $\Uc$ $\in$ $\Sc_{d,\ell,m,k,d'}^{\mfs,N,\rho}$,  is an $\R^{d'}$-valued exchangeable function to the order $N$ on $\R^d$,  in the form: 
\begin{align} \label{defsymNN}
\Uc(\bolx) &= \; \psi\big( \mfs((\varphi(x_i))_{i\in \llbracket 1,N\rrbracket}) \big),  \quad  \mbox{ for } \bolx = (x_i)_{i \in \llbracket 1,N\rrbracket} \in (\R^d)^{N},    
\end{align}
where $\varphi$ $\in$ $\Nc_{d,\ell,m,k}^\rho$, $\psi$ $\in$ $\Nc_{k,\ell,m,d'}^\rho$ (here, for simplicity of notations, we assume that the number of hidden layers and neurons  of $\varphi$ and $\psi$ are the same but in practical implementation,  they may be different), 
and $\mfs$ is a given $\R^k$-valued exchangeable function to the order $N$ on $\R^k$, typically: 
\begin{itemize}
\item Max-pooling (component-wise):  $\mfs(\boly)$ $=$ $\max(y_i)_{i \in \llbracket 1,N\rrbracket}$,  
\item Sum: $\mfs(\boly)$ $=$ $\sum_{i=1}^N y_i$, or average:  $\mfs(\boly)$ $=$ $\frac{1}{N}\sum_{i=1}^N y_i$,
\end{itemize}
for $\boly$ $=$ $(y_i)_{i\in \llbracket 1,N \rrbracket}$ $\in$ $(\R^k)^N$. 
When $\mfs$ is the max-pooling function, $\Sc_{d,\ell,m,k,d'}^{\mfs,N,\rho}$ is called {\it PointNet}, as introduced in \cite{pointnet}, while for $\mfs$  equals to the sum/average function,  it is called {\it DeepSet}, see \cite{deepsets}.  The 
architecture is described in Figure \ref{fig: symNN}, and $k$ can be interpreted as a number of features describing the geometry of the set of points $\{x_i\}_{i\in\llbracket 1,N\rrbracket}$. For example in the context of mean-field control problem, $k$ will be related to the moments for describing the 
law of the McKean-Vlasov SDE.

 A given symmetric network  function $\Uc$ $\in$ $\Sc^{\mfs,N,\rho}_{d,\ell,m,k,d'}$  is determined by the weight/bias parameters $\theta$ $=$ $(\theta^{(1)},\theta^{(2)})$ with 
$\theta^{(1)}$ $=$ $(\Wc_0^{(1)},\beta_0^{(1)},\ldots,\Wc_\ell^{(1)},\beta_\ell^{(1)})$ defining the layer functions $\phi_0 \ldots,\phi_\ell$ of $\varphi$, 
and $\theta^{(2)}$ $=$ $(\Wc_0^{(2)},\beta_0^{(2)},\ldots,\Wc_\ell^{(2)},\beta_\ell^{(2)})$ defining the layer functions $\psi_0 \ldots,\psi_\ell$ of $\psi$. The number of parameters is 
$M$ $=$ $M_1+M_2$, with  $M_1$ $=$ $m(d+1) +  m(m+1)(\ell-1) + (m+1)k$, $M_2$ $=$ $(k+1)m +  m(m+1)(\ell-1) + (m+1)d'$, and we observe that it does not depend on  
the number $N$ of inputs.  

\begin{Remark}[{\it Time dependent symmetric network}]
{\rm 
A  time-dependent symmetric in space neural network can be constructed as
\begin{align}
\Uc(t,\bolx) &=&  \psi\big(t,\mfs((\varphi(x_i))_{i\in \llbracket 1,N\rrbracket}) \big),  \quad  \mbox{ for } t \in \R_+,  \;  \bolx = (x_i)_{i \in I} \in (\R^d)^{N},
\end{align}
with $\varphi$ a feedforward network from $\R^d$ into $\R^k$, and $\psi$ is a feedforward from $\R^{k+1}$ into $\R^{d'}$, where we add time as an additional feature, see architecture in Figure \ref{fig:timeNN}. 

 \begin{figure}[H]
    \centering
    \includegraphics[width=10cm]{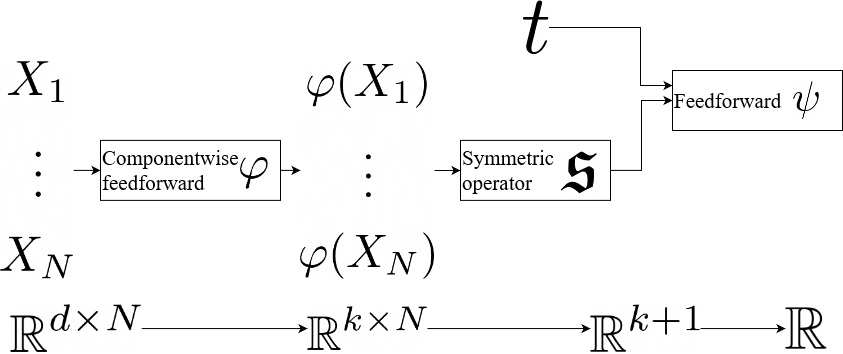}
    \caption{Architecture  of time dependent symmetric network.
    }
    \label{fig:timeNN}
\end{figure}
}
\ep
\end{Remark}

\begin{figure}[H]
    \centering
    \includegraphics[width=13cm]{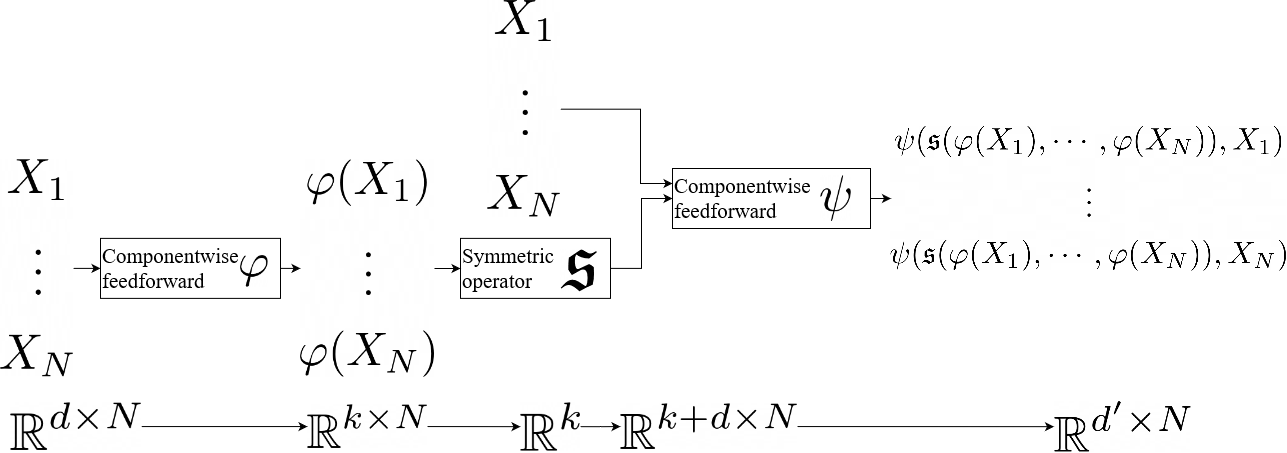}
    \caption{Architecture of DeepDerSet network}
    \label{fig:DeepDerSet}
\end{figure}

A {\it $D$-symmetric neural  network} function, denoted $\Zc$ $\in$ $D\Sc^{\mfs,N,\rho}_{d,\ell,m,k,d'}$,  
is an $\R^{d'}$-valued  $D$-exchan\-geable function in the form
\begin{align} \label{defZc} 
\Zc(\bolx,x) &= \; \psi\big( \mfs((\varphi(x_i))_{i\in  \llbracket 1,N\rrbracket}),x \big),  \quad \mbox{ for } \bolx = (x_i)_{i \in  \llbracket 1,N\rrbracket} \in (\R^d)^{N}, x \in \R^d,   
\end{align}
where $\varphi$ $\in$ $\Nc_{d,\ell,m,k}^\rho$, $\psi$ $\in$ $\Nc_{k+d,\ell,m,d'}^\rho$, and $\mfs$ is a given $\R^k$-valued $N$-exchangeable function on $\R^k$. 
The number of parameters of a given $\Zc$ $\in$ $D\Sc^{\mfs,\rho}_{d,\ell,m,k,d'}$ is $M'$ $=$ $M'_1+M'_2$, with 
$M'_1$ $=$ $m(d+1) +  m(m+1)(\ell-1) + (m+1)k$, $M'_2$ $=$ $m(k+d+1) +  m(m+1)(\ell-1) + (m+1)d'$. 
We shall often take for $\mfs$ the average function, and call $D\Sc^{\mfs,N,\rho}_{d,\ell,m,k,d'}$ as {\it DeepDerSet}. Its architecture is given in Figure \ref{fig:DeepDerSet}. 
Given a $D$-symmetric neural network $\Zc$, we denote by $\boZc$ the function defined on $(\R^d)^N$ by $\boZc(\bolx)$ $=$ $(\Zc(\bolx,x_i))_{\in\llbracket 1,N\rrbracket}$ $\in$ $(\R^d)^N$, and by misuse of notation, we may also call $\boZc$ as a $D$-symmetric NN. By construction, these networks respect the representation given by Lemma \ref{lemvNderiv}, by being defined as a $D$-exchangeable function applied component by component. In that way we are able to enforce the correct symmetries for representing both a symmetric function and its derivative, which will be useful in Section  \ref{secalgo} and Section \ref{secnum} when looking for the gradient of PDE solutions. 
  
 \vspace{1mm}

 Alternatively,  one can generate $D$-exchangeable functions as follows. Starting from a DeepSet element $\Uc$ $\in$  $\Sc_{d,\ell,m,k,d'}^{\mfs,N,\rho}$ as in \eqref{defsymNN} with $\mfs$ the sum function, and network functions $\varphi,\psi$ with smooth activation functions
\begin{align}
D_{x_i} \Uc(\bolx) &=  {\rm D}\Uc(\bolx,x_i), \quad \bolx = (x_i)_{\in \llbracket 1,N\rrbracket} \in (\R^d)^N,  
\end{align}  
with 
\begin{align} \label{AD-DeepSet}
{\rm D}\Uc(\bolx,x) & := \;  D\varphi(x) D\psi\big( \mfs((\varphi(x_i))_{i\in I}) \big), \quad  x \in \R^d. 
\end{align}  
The set of $D$-exchangeable functions obtained from  differentiation of DeepSet network functions is called  {\it AD-DeepSet}, where AD stands for automatic differentiation, and Automatic refers to the implementation of the differentiation in software library, e.g. in TensorFlow. This alternative network for approximating the derivative of a differentiable symmetric function naturally respects the structure given by Lemma \ref{lemvNderiv}, by construction as the derivative of a symmetric DeepSet neural network.

Given a $D$-symmetric NN $\boZc$ in {\it DeepDerSet} or in {\it AD-DeepSet} with smooth activation functions, we denote by $D\boZc$ its  differentiation. 

\vspace{1mm}

 As for the well-known universal approximation theorem \cite{hor91} for neural networks, we have a similar result for symmetric neural networks, which states that  any  exchangeable function can be arbitrarily approximated by a PointNet or DeepSet  given enough neurons. 
 More precisely, by combining Theorem 2.9 of \cite{limitations} with Theorem 2 of \cite{hor91}, we obtain the following approximation theorem for DeepSets.
 
 \vspace{2mm}
 
 \noindent {\bf Universal approximation for DeepSets networks.} 
 Let $\mfs$ be the sum function.  The set $\cup_{m=1}^\infty \Sc^{\mfs,N,\rho}_{d,\ell,m,N+1,d'}$  
 approximates any $N$-exchangeable continuous function $w$ 
arbitrary well on any compact  set of  $K$ $\subset$ $\R^d$, once $\rho$  is continuous, bounded and non-constant:  for all  $\eps$ $>$ $0$, $N$ $\in$ $\N^*$, there exists $\Uc$ $\in$ $\cup_{m=1}^\infty \Sc^{\mfs,N,\rho}_{d,\ell,m,N+1,d'}$ such that
\begin{align}
 \big| w(\bolx) - \Uc(\bolx) \big| & \leq \; \eps \quad \forall  \bolx  \in K^N.
\end{align}
Note that a priori the latent space dimension $k$ has to be chosen equal to  $N+1$.

\vspace{1mm}

Alternatively, by combining Theorem 1 of \cite{pointnet} with Theorem 2 of \cite{hor91},  we obtain the following one-dimensional approximation theorem for PointNet. 

 \vspace{2mm}

\noindent {\bf Universal approximation for PointNet networks.} 
Let $\mfs$ be the $\max$ function. The set $\cup_{m=1}^\infty \cup_{k=1}^\infty \Sc^{\mfs,N,\rho}_{1,\ell,m,k,d'}$  
 approximates any $N$-exchangeable Hausdorff continuous function $w$ (seen as a function on sets)
arbitrary well on any compact  set of  $K$ $\subset$ $\R$
, once $\rho$  is continuous, bounded, and non-constant:  for all $\eps$ $>$ $0$, $N$ $\in$ $\N^*$, there exists $\Uc$ $\in$ $\cup_{m=1}^\infty \cup_{k=1}^\infty \Sc^{\mfs,N,\rho}_{1,\ell,m,k,d'}$ such that
\begin{align}
\big| w(S) - \Uc(\bolx) \big| & \leq \; \eps, \quad \forall S\subset K,\ S=\{x_1,\cdots,x_N\}. 
\end{align}
Note here that a priori the latent space dimension $k$ has to be chosen as large as needed.

\subsection{Comparison tests}

In this paragraph, we test the accuracy of the approximation of exchangeable functions by {\it DeepSet} or {\it PointNet},  and also the approximation of $D$-exchangeable functions by {\it DeepDerSet} or {\it  AD-DeepSet}, and compare numerically with 
classical feedforward approximations.

 \subsubsection{Approximation of some simple functions} \label{sec:approxEx}

We first test  the approximation of the following simple symmetric functions:
 \begin{enumerate}
     \item  $f(x)$ $=$ $\exp{(2 \bar x + 3 \bar{x}^3)}$,  with  $\bar x$ $=$ $\displaystyle{\frac{1}{N} \sum_{i=1}^N  x_i}$    (case 1)
     \item $f(x)$ $=$   $\displaystyle{\frac{1}{N} \sum_{i=1}^N} \big[  \sin(x_i)  1_{x_i <0} + x_i 1_{x_i \ge 0}\big]$    (case 2)
     \item $f(x)$ $=$  $\bar x + 2 \bar{x}^2 + 3 \bar{x}^3$,  with  $\bar x$ $=$ $\max\{ x_i, i=1, \dots, N\}$   (case 3)
     \item $f(x)$ $=$ $\cos( 2\bar{x} + 3 \bar{x}^2)$, with  $\bar x$ $=$ $ \sum_{i=1}^N x_i$ (case 4)
 \end{enumerate}
 
We  use a symmetric neural network architecture as proposed in  \cite{pointnet},  \cite{deepsets}:
 \begin{itemize}
 \item First, 
a feedforward network $\varphi$ with $\ell$ $=$ $5$ hidden layers,  and respectively $64,64,64,128$ and $1024$ neurons such that each dimension $i$, $i=1,\dots, N$,  is treated with the same network in one dimension avoiding to break the symmetry. 
 \item Two possible symmetric functions $\mfs$ to the order $N$ on $\R^k$ with $k$ $=$ $1024$, the max-pooling ({\it PointNet}) and the sum function ({\it DeepSet}).  
 \item At last, a feedforward network $\psi$ from $\R^{1024}$ to $\R$ with $\ell$ $=$ $2$  hidden layers, and respectively $512$ and $256$ neurons. 
 \end{itemize}
 
 For the approximation with classical feedforward networks, we used three or four layers and a number of neurons constant per layer  equal to $10 +d$, or $10+2 d$ neurons.

We use the ADAM  optimizer (\cite{kingma2014adam}), with a batch size equal to  $300$  for solving the approximation problem with quadratic loss function:  
 \begin{flalign}
  \min_{\theta} \E[ |f(X) - \Uc^{\theta}(X)|^2],
  \label{eq:error space}
 \end{flalign}
 with training simulations from  $X$ $\sim$ $\Nc(0_N,1_N)$, and  $\theta$ are the  parameters of the network function $\Uc^\theta$. 
 The number of epochs (corresponding to the number of gradient descent iterations) used  is equal to $100$.  After epoch iterations of the
 stochastic gradient, the error \eqref{eq:error space} is estimated with $20000$ simulations. If the error is below a threshold equal to \num{1e-5} the optimization is stopped, otherwise a counter for outer iterations is incremented.
 The number  of outer iterations is blocked at epochExt $=1000$ (meaning a maximal total number of stochastic gradient iterations equal to epoch $\times$ epochExt $= 100000$).
 
 In Tables \ref{tab:func1}, \ref{tab:func2}, \ref{tab:func3}, we report the accuracy reached (Error) and  the number of iterations (Iter.)  used to obtain this given accuracy: then  a threshold equal to  \num{1e-5} means that the optimization has been successful and the relevant parameter is the number of iterations used. A number of iterations equal to expochExt $=$ $1000$  means that the optimization has not been successful and the error reached indicates how far we are from optimality.  For the feedforward case, we report the best result (``minimum" in table) and the worse result (``maximum" in table) obtained changing the number of layers and the number of neurons used.

 The initial learning rate is taken equal to \num{1e-3} for first outer simulation in cases 1 and 2 with a linear decay to  \num{1e-5} for a number of outer iterations equal to $1000$. For test case 3, the initial learning rate is taken equal to \num{1e-4} with a linear decay to  \num{1e-5}.
 The result obtained in Table  \ref{tab:func1} is similarly obtained with a large number of functions tested in dimension between $N$ $=$ $10$ to $1000$. 
 It shows the following results:
 \begin{itemize}
     \item The classical feedforward, with dense layers, often permits to obtain optimally without forcing symmetry of the solution,
     \item Classical feedforward results do not depend a lot on the number of layers, the number of  neurons tested and the activation function used,
     \item For symmetric approximations,  DeepSets generally permits to get the best results and the ReLU activation  function is the best out of the three tested.
 \end{itemize}
 \begin{table}[h!]
     \centering
     \begin{tabular}{|c|c|c|c|c|c|c|c|c|} \hline
         &  \multicolumn{4}{c|}{Symmetric}  &  \multicolumn{4}{c|}{ Feedforward} \\ \hline
       & \multicolumn{2}{c|}{PointNet}  & \multicolumn{2}{c|}{DeepSet}  &  \multicolumn{2}{c|}{Minimum} &  \multicolumn{2}{c|}{Maximum}  \\ \hline
        Activation  &  Error & Iter. & Error & Iter. & Error & Iter. & Error  & Iter. \\ \hline
         ReLU          &  $0.008$  & $1000$ &  \num{1e-5}  & $10$ &  \num{1e-5}   & $125$  &  \num{1e-5} & $166$  \\ \hline
         tanh         & $0.016$   & $1000$  &  \num{1e-5}   & $288$  &  \num{1e-5}   & $180$  &  \num{1e-5} & $308$\\ \hline
         ELU          & $0.015$   & $1000$  &  \num{1e-5}   & $176$  &  \num{1e-5}   & $108$  &  \num{1e-5} & $130$\\ \hline  
    \end{tabular}
     \caption{Approximation error \eqref{eq:error space} obtained for different networks on one run and number of iterations used depending on  activation functions for approximation of the function $f$ in case 1, dimension $N$ $=$ $100$.}
     \label{tab:func1}
 \end{table}
 
In the sequel, we drop the ELU activation function on other cases as shown   in Tables \ref{tab:func2} for cases 2 and 3.  Notice that case 3, involving a max function is the only one where PointNet approximation gives the best result among the other tested  networks.
On cases 2 and 3 in dimension 100,  the DeepSets approximation  outperforms the classical feedforward network for all the number of layers and neurons tested. However, results on case 3 are not very good for the PointNet approximation even with the ReLU activation function.\\
 \begin{table}[h!]
     \centering
         \begin{tabular}{|c|c|c|c|c|c|c|c|c|c|} \hline
       Case &  activ & \multicolumn{4}{c|}{Symmetric}    &  \multicolumn{4}{c|}{ Feedforward} \\ \hline
     &   & \multicolumn{2}{c|}{PointNet}  & \multicolumn{2}{c|}{DeepSets}  &  \multicolumn{2}{c|}{Minimum} &  \multicolumn{2}{c|}{Maximum}  \\ \hline
          &  &  Error & Iter. & Error & Iter. & Error & Iter. & Error & Iter. \\ \hline
        2    & ReLU      &  $0.001$  & $1000$ &  \num{1e-5}  & $5$ &  \num{1e-5}   & $992$  & $0.002$ & $1000$  \\ \hline
         2     & tanh     &  $0.03$  & $1000$ &  \num{1e-5}  & $342$ & $0.0018$   & $1000$  & $6\times$  \num{1e-5} & $1000$  \\ \hline
        3   &  ReLU      &  $0.001$  & $1000$ & $0.23$  & $1000$ & $88$   & $1000$  & $432$ & $1000$  \\ \hline
          3   &  tanh     &  $0.002$  & $1000$ & $65$  & $1000$ & $933$   & $1000$  & $969$ & $1000$  \\ \hline
    \end{tabular}
     \caption{Approximation error \eqref{eq:error space} obtained for different networks on one run and number of iterations used  for approximation of the function $f$ in cases 2 and 3, dimension $N$ $=$ $100$, activation function ReLU.}
     \label{tab:func2}
 \end{table}
 
Results for test case 4 are given on Table \ref{tab:func3} using an initial learning rate equal to \num{5e-5} and a decay linear to \num{5e-6} with the number of outer iterations.
 \begin{table}[h!]
  \centering
         \begin{tabular}{|c|c|c|c|c|c|c|c|c|} \hline
        &\multicolumn{4}{|c|}{Symmetric}  &  \multicolumn{4}{c|}{ Feedforward} \\ \hline
       & \multicolumn{2}{|c|}{PointNet}  & \multicolumn{2}{c|}{DeepSets}  &  \multicolumn{2}{c|}{Minimum} &  \multicolumn{2}{c|}{Maximum}  \\ \hline
        activ  & Error & Iter. & Error & Iter. & Error & Iter. & Error  & Iter. \\ \hline
       ReLU  &  $0.1910$  & $1000$ & \num{8.9e-5}  & $1000$ & $0.0045$   & $1000$  & $0.01$ & $1000$  \\ \hline
              tanh &  $0.19$  & $1000$ & \num{4e-5}  & $1000$ & $0.0006$   & $1000$  & $0.0012$ & $1000$  \\ \hline
    \end{tabular}
     \caption{Approximation error \eqref{eq:error space} obtained for different networks, activation functions for approximation of the function $f$ in case 4 dimension 1000.}
     \label{tab:func3}
 \end{table}
 At last, considering case 4 in dimension $N$ $=$ $1000$, when the function is quickly changing, we see that the classical feedforward network functions have difficulty to converge while the DeepSets network  approximation converges. The latter turns out to 
be a very good candidate to some very high dimensional PDEs when there is symmetry in the solution.

 \subsubsection{Approximation of a function of $t$ and $x$ with symmetry in $x$}  \label{sec:tAndX}

We test the accuracy of our time dependent symmetric neural network by considering the following two cases of functions:  
\begin{enumerate}
     \item $f(x)= \exp\big(\bar{x} (t + 2t^2) + 3 t \bar{x}^3\big) $  with $\bar x$ $=$ $\displaystyle{\frac{1}{N} \sum_{i=1}^N x_i}$   (case 1)
     \item $f(x)= t+ \cos( t \bar{x})$, with $\bar x$ $=$ $\frac{1}{\sqrt{N}} \sum_{i=1}^N x_i$  (case 2)
 \end{enumerate}
The approximation is performed through the minimization problem
 \begin{flalign}
  \min_{\theta} \E[ |f(\tau,X) - \Uc^{\theta}(\tau,X)|^2],
  \label{eq:error time space}
 \end{flalign}
 with training simulations from  $X \sim \Nc(0_N,1_N)$, and an independent uniform law for $\tau$ on $[0,1]$, and where $\Uc^\theta$ is a time-dependent DeepSet with parameters $\theta$. We keep the same number of neurons and layers as in the previous section, and compare  
with  a classical feedforward network  composed of $3$ layers of  $d+10$ neurons.  In all experiments,  we use a ReLU activation function.

 \begin{table}[h!]
     \centering
         \begin{tabular}{|c|c|c|c|c|} \hline
       Case &  \multicolumn{2} {c|}{DeepSets}  &  \multicolumn{2}{c|}{ Feedforward} \\ \hline
          &  Error & Iter. & Error & Iter.  \\ \hline
         1  &   \num{1e-5}  & $67$ & $0.008$  & $1000$  \\ \hline
         2   &   \num{1e-5}  & $344$ & $0.048$  & $1000$   \\ \hline
    \end{tabular}
     \caption{Approximation error \eqref{eq:error time space} obtained for different networks on one run and number of iterations used  for approximation of the function $f$ in case 1 and 2, dimension $N$ $=$ $100$.}
     \label{tab:funcWithT1}
 \end{table}

In Table \ref{tab:funcWithT1}, we give the results obtained in dimension 100. Surprisingly, the feedforward approximation seems to have difficulties to approximate the case 1 although it is quite similar to case one in the previous section.
 As for the second case, the result is not so surprising as the case is quite similar to case 4 in previous section, where the feedforward network has difficulties to converge.

 \subsubsection{Gradient approximation}

 We now focus on  the approximation of the derivative of an exchangeable function by means of $\Uc^{\theta}$ a {\it DeepDerSet}, a {\it  AD-DeepSet}, or a classical feedforward network . 

 The minimization problem is now:
 \begin{flalign}
  \min_{\theta} \E[ \|Df(X) - \Uc^{\theta}(X)\|^2],
  \label{eq:error Grad}
 \end{flalign}
 where  the norm $\| \cdot \|$ is  the Euclidean norm.

 The comparison is performed on the following test functions: 

 \begin{enumerate}
     \item  $f(x)$ $=$  $\exp{( \bar x + \bar{x}^3)} ( 1_N + 3 x^2) $ where $\bar x$ $=$ $\displaystyle{\frac{1}{N} \sum_{i=1}^N x_i} $   (case 1)
     \item $f(x)$ $=$ $y$ where $y_i=  1_{x_i >0} + cos(x_i) 1_{x_i< 0}, \quad i=1, \dots, N$   (case 2)
     \item $f(x)$ $=$  $\frac{1}{\sqrt{N}} \sin( \bar{x}) 1_N$, where $\bar x$ $=$ $\frac{1}{\sqrt{N}}  \sum_{i=1}^N x_i$ (case 3)
 \end{enumerate}
 
 We compare the classical feedforward approximation to our network approximation in Tables   \ref{tab:funcWithGradRelu}  and \ref{tab:funcWithGradTanh}, using a maximal number of iterations equal to $5000$.  
 Clearly  using a ReLU activation function is superior to the tanh activation function and  the {\it DeepDerSet}  gives the best approximation while the  {\it AD-DeepSet}  or the feedforward may have difficulties  to approximate the functions accurately.
\begin{table}[h!]
     \centering
         \begin{tabular}{|c|c|c|c|c|c|c|c|} \hline
       Case & $N$ & \multicolumn{2} {c|}{Feedforward}  &  \multicolumn{2}{c|}{DeepDerSet }& \multicolumn{2}{c|}{AD-DeepSet } \\ \hline
          &   & Error & Iter. & Error & Iter. & Error  & Iter.  \\ \hline
         1  &  10 & \num{5e-4}  & 5000 &   \num{5e-4}  & 5000 &   \num{7e-3}  & 5000  \\ \hline
         1  &  100 & \num{1e-5}  & 300 & \num{1e-5} &  250 & \num{1e-5} &  50  \\ \hline
         2   & 10 &  \num{0.03}  & 5000 &  \num{1e-5}  &  550 &  \num{2e-4}  & 5000\\ \hline
         2   &  100 &  \num{0.12}  &  5000  & \num{1e-5}  &  450 &  \num{1e-4}  & 5000\\ \hline
         3   & 10 &  \num{1e-5}  & 3800 &  \num{1e-5} &  850 & \num{0.03}  & 5000 \\ \hline
          3   &  100 &  \num{1e-5}  & 1850 &  \num{1e-5} &  700 & \num{3E-3}  & 5000 \\ \hline
  \end{tabular}
     \caption{Approximation error \eqref{eq:error Grad} with ReLU activation function obtained for different networks on one run and number of iterations used  for approximation of the derivative of an exchangeable function.}
     \label{tab:funcWithGradRelu}
 \end{table}
 
\begin{table}[h!]
     \centering
         \begin{tabular}{|c|c|c|c|c|c|c|c|} \hline
       Case & $N$ & \multicolumn{2} {c|}{Feedforward}  &  \multicolumn{2}{c|}{DeepDerSet }& \multicolumn{2}{c|}{AD-DeepSet } \\ \hline
          &   & Error & Iter. & Error & Iter. & Error  & Iter.  \\ \hline
         1  &  10 & \num{5e-4}  & 5000 &   \num{1.7e-4}  & 5000 &   \num{3e-3}  & 5000  \\ \hline
         1  &  100 & \num{1e-5}  & 500 & \num{1e-5} &  100 & \num{1e-5} &  50  \\ \hline
         2   & 10 &  \num{0.03}  & 5000 &  \num{1e-5}  &  950 &  \num{0.028}  & 5000\\ \hline
         2   &  100 &  \num{0.12}  &  5000  & \num{1e-5}  &  700 &  \num{0.56}  & 5000\\ \hline
         3   & 10 &  \num{1.5e-5}  & 5000 &  \num{1e-5} &  1850 & \num{0.02}  & 5000 \\ \hline
          3   &  100 &  \num{1e-5}  & 2100 &  \num{1e-5} &  1350 & \num{4E-3}  & 5000 \\ \hline
  \end{tabular}
     \caption{Approximation error \eqref{eq:error Grad} with tanh  activation function obtained for different networks on one run and number of iterations used  for approximation of the derivative of an exchangeable function.}
     \label{tab:funcWithGradTanh}
 \end{table}

\section{Numerical schemes}  \label{secalgo}

We now adapt the deep backward dynamic programming (DBDP) schemes  developed in \cite{hure2020deep}  and \cite{phawarger20} for solving nonlinear PDEs, by using symmetric neural networks and $D$-symmetric neural networks (instead of feedforward neural networks) for approximating 
the exchangeable solution $v$ and  its gradient $D_{\bolx} v$.   We recall the main steps of the DBDP scheme, and distinguish  the case of semi-linear and fully non-linear PDEs.

\subsection{Semi-linear PDE}

We first consider the  case where the generator $F$ in \eqref{symPDE}  may be  decomposed into the form 
\beq \label{semiF}
F(t,\bolx,y,\boz,\bog) &=& H(t,\bolx,y,\boz) + \sum_{i=1}^N b_i(t,\bolx).z_i +  \frac{1}{2} \sum_{i,j=1}^N {\rm tr}\big(\Sigma_{ij}(t,\bolx) \gamma_{ij} \big), 
\enq
for $t$ $\in$ $[0,T]$, $\bolx$ $=$ $(x_i)_{i\in\llbracket 1,N\rrbracket}$ $\in$ $(\R^d)^N$, $y$ $\in$ $\R$, $\boz$ $=$ $(z_i)_{i\in\llbracket 1,N\rrbracket}$ $\in$ $(\R^d)^N$, and 
$\bog$ $=$ $(\gamma_{ij})_{i,j\in\llbracket 1,N\rrbracket}$ $\in$ $\S^N(\S^d)$. 
Here,  $H$ is a function on $[0,T]\times(\R^d)^N\times\R\times(\R^d)^N$ satisfying the permutation-invariance condition:  
\beqs
H(t,\bolx,y,\boz) &=& H(t,\pi[\bolx],y,\pi[\boz]),   \quad \forall \pi \in  \mfS_N,
\enqs
the coefficients $b_i$, $i$ $=$ $1,\ldots,N$, are $\R^d$-valued functions on $[0,T]\times(\R^d)^N$ satisfying the condition
\begin{align} \label{symb}
b_i(t,\pi[\bolx]) &=\;  b_{\pi(i)}(t,\bolx),  \quad \forall \pi \in  \mfS_N,
\end{align}
and the coefficients $\Sigma_{ij}$, $i,j$ $=$ $1,\ldots,N$,  are $d\times d$-matrix valued functions on $[0,T]\times(\R^d)^N$ in the form 
\begin{align} \label{defSig} 
\Sigma_{ij}(t,\bolx) &= \;  \sigma_{ij}\sigma_{ij}\trans(t,\bolx)  +  \sigma_{i0}\sigma\trans_{j0}(t,\bolx), 
\end{align}
 for some $d\times d$-matrix valued functions $\sigma_{ij}$, and $d\times q$-matrix valued functions $\sigma_{i0}$, satisfying the invariance property: for all $\pi$ $\in$ $\mfS_N$,
 \begin{equation} \label{symsigbis} 
 \sigma_{ij}(t,\pi[\bolx]) \; = \;  \sigma_{\pi(i)\pi(j)}(t,\bolx), \quad \quad  \sigma_{i0}(t,\pi[\bolx]) \; = \;  \sigma_{\pi(i)0}(t,\bolx). 
 \end{equation}
In this case, the permutation-invariance condition {\bf (HI)} on $F$ is  satisfied, and we observe that it includes Example \ref{exmultiasset} of multi-asset pricing with 
$H(t,\bolx,y,\boz)$ $=$ $\beta(y^+-y)$ (in the case of the CVA pricing), $b_i$ $\equiv$ $r$  and $\sigma_{i0}$ $\equiv$ $0$. This also 
includes Example \ref{exMKV} of the McKean-Vlasov control problem under common noise with uncontrolled diffusion coefficient, where $\sigma_{ij}(t,\bolx)$ $=$ $\sigma(t,x_i,\bar\mu(\bolx)) \delta_{ij}$, 
$\sigma_{i0}(t,\bolx)$ $=$ $\sigma_0(t,x_i,\bar\mu(\bolx))$,  and   
\begin{align} \label{HMKV}
H(t,\bolx,y,\boz) &=   -ry+ \sum_{i=1}^N \inf_{a \in A} \big[  \beta(t,x_i,\bar\mu(\bolx),a).z_i +  \frac{1}{N}f(x_i,\bar\mu(\bolx),a) \big] \\
& \quad \quad \quad - \sum_{i=1}^N  b_i(t,\bolx).z_i,   
\end{align}
for any function $b_i$  satisfying \eqref{symb}.  
We shall discuss more in detail    the relevant choice of the drift coefficient $b_i$ in Section \ref{sec: extension}. 

\vspace{2mm}

The starting point of the numerical scheme is the probabilistic representation of the PDE \eqref{symPDE} with $F$ as in \eqref{semiF} in terms of a forward backward stochastic differential equation (FBSDE), as in \cite{pardoux1990adapted}.   
In our context, the forward system is described by the process $\boX$ $=$ $(X^1,\ldots,X^N)$ valued in $(\R^d)^N$ governed by the diffusion dynamics: 
\begin{align} \label{diffX}
dX_t^i &=\;  b_i(t,\boX_t) dt  + \sum_{j=0}^N \sigma_{ij}(t,\boX_t) dW_t^j,  
\end{align}
where $W^i$, $i$ $=$ $1,\ldots,N$, are independent $d$-dimensional Brownian motions, independent of the  $q$-dimensional Brownian motion $W^0$.  Given this forward diffusion process, we then consider the pair process 
$(Y,\boZ=(Z^i)_{i\in \llbracket 1,N\rrbracket})$  valued in $\R\times(\R^d)^N$ solution to the BSDE
\begin{align} \label{BSDEsemi}
G(\boX_T) - Y_t +  \int_t^T H(s,\boX_s,Y_s,\boZ_s) ds &  \\
-  \sum_{i=1}^N \sum_{j=0}^N  \int_t^T  (Z_s^i)\trans  \sigma_{ij}(s,\boX_s) dW_s^j  & \; =\; 0, \quad \quad 0 \leq t\leq T,  
\end{align}
which is connected by It\^o's formula  to the solution of the PDE  \eqref{symPDE} via: 
\beqs
Y_t \; = \; v(t,\boX_t), & & Z_t^i \; = \; D_{x_i} v(t,\boX_t), \; i=1,\ldots,N, \quad 0 \leq t \leq T. 
\enqs

We next consider a time discretization of this FBSDE on a time grid $\{t_k, k=0,\ldots,N_T\}$, with $t_0$ $=$ $0$, $t_{N_T}$ $=$ $T$, 
$\Delta t_k$ $:=$ $t_{k+1}-t_k$ $>$ $0$,  by defining the Euler scheme $\{\boX_k^{N_T} = (X_k^{i,N_T})_{i\in\llbracket 1,N\rrbracket}$, $k=0,\ldots,N_T\}$ associated to  the forward diffusion process  
$\{\boX_t=(X_t^i)_{i\in\llbracket 1,N\rrbracket}$, $0\leq t \leq T\}$, which is used for the training simulations, together with the increments of the Brownian motions:  $\Delta W_k^j$ $:=$ $W_{t_{k+1}}^j-W_{t_k}^j$, $k$ $=$ $0,\ldots,N_T-1$, $j$ $=$ $0,\ldots,N$,   
of our numerical backward scheme.   The DBDP algorithm reads then as follows:

\begin{algorithm2e}[h!] 
\DontPrintSemicolon 
\SetAlgoLined 
\vspace{1mm}
{\bf Initialization}: Initialize from the exchangeable function: $\widehat\Uc_{N_T}(\cdot)$ $=$ $G(\cdot)$ 

\For{ $k$ $=$ $N_T-1,\ldots,0$} 
{minimize over  symmetric NN  $\Uc_k$, and $D$-symmetric NN  $\Zc_k$, the quadratic loss function
\begin{align} \label{loss}
J_1(\Uc_k,\Zc_k) &= \; \E\Big| \widehat\Uc_{k+1}(\boX_{k+1}^{N_T}) - \Uc_k(\boX_k^{N_T})  \\
&  \quad \quad \quad  + \;   H\big(t_k,\boX_k^{N_T},\Uc_k(\boX_k^{N_T}),\boZc_k(\boX_k^{N_T})\big) \Delta t_k  \\
& \quad \quad  \quad - \;  \sum_{i=1}^N \sum_{j=0}^N \big(\Zc_k(\boX_k^{N_T},X_k^{i,N_T})\big) \trans\sigma_{ij}\big(t_k,\boX_k^{N_T}\big) \Delta W_k^j  \Big|^2,
\end{align}
and update $(\widehat\Uc_k,\widehat\Zc_k)$ as the solution to this minimization problem. 
}
\caption{DBDP scheme with symmetric NN  \label{MKVscheme} }
\end{algorithm2e}

The output of the DBDP scheme provides  approximations $\widehat\Uc_k(\bolx)$ of $v(t_k,\bolx)$,  and $\widehat\boZc_k(\bolx)$ of $D_{\bolx} v(t_k,\bolx)$, $k$ $=$ $0,\ldots,N_T-1$, for values of $\bolx$ $\in$ $(\R^d)^N$ that are 
well-explored by the training simulations of $\boX_k^{N_T}$.   We refer to \cite{phawarger20} (see their section 3.1) for a discussion on the choice of the algorithm hyperparameters. 

\begin{Remark}
We stress that the neural networks do not take time as an input. Adding time would not make any difference because the training is done locally, time step per time step. So the neural networks approximating $\Uc_k$ and $\Zc_k$ would be trained with only $t_{k}$ as an input and hence they would not be able to learn the dependence on time. However, at time $k < N_T$, we initialize the parameters of $\Uc_k$ and $\Zc_k$ respectively with the parameters of the neural networks for $\Uc_{k+1}$ and $\Zc_{k+1}$, which have been trained at the previous iteration. This gives a good initial guess for the neural networks at time $t_{k}$ and leads to more efficient training. 
\end{Remark}
 
\subsection{Fully nonlinear PDE}  \label{secfull}

We consider more generally the fully non-linear PDE \eqref{symPDE} with a symmetric generator $F$ satisfying {\bf (HI)}. 
We  adapt  the machine learning scheme in \cite{phawarger20} for solving fully nonlinear PDEs by exploiting furthermore the exchangeability property of the solution by using again symmetric neural networks as in the semi-linear case. 

We fix some arbitrary drift and diffusion coefficients $b_i$, $\sigma_{ij}$, $i$ $=$ $1,\ldots,N$, $j$ $=$ $0,\ldots,N$, satisfying invariance properties as in \eqref{symb}-\eqref{symsigbis}  
(in practice, they should be chosen depending on the studied problem as for the semi-linear case, see a general discussion in  Section 3.1 in \cite{phawarger20}, and an application in Section \ref{sec:MeanVar}),   
and introduce the forward diffusion system $\boX$ as in \eqref{diffX} and  its discrete-time Euler scheme $\boX^{N_T}$.    
We then consider the triple process $(Y,\boZ=(Z^{i})_{i\in\llbracket 1,N\rrbracket},\boG=(\Gamma^{ij})_{i,j\in\llbracket 1,N\rrbracket})$ valued in $\R\times(\R^d)^N\times\S^N(\S^{d})$ solution to the BSDE
\begin{align}
G(\boX_T) - Y_t +  \int_t^T F_{b,\sigma}(s,\boX_s,Y_s,\boZ_s,\boG_s) ds  & \\
-  \sum_{i=1}^N \sum_{j=0}^N  \int_t^T  (Z_s^i)\trans  \sigma_{ij}(s,\boX_s) dW_s^j  &=\; 0, \quad 0 \leq t\leq T,  
\end{align}
with  
\begin{align}
F_{b,\sigma}(t,\bolx,y,\boz,\bog) &:= \;  F(t,\bolx,y,\boz,\bog)  - \sum_{i=1}^N b_i(t,\bolx).z_i 
- \frac{1}{2} \sum_{i,j=1}^N {\rm tr}\big(\Sigma_{ij}(t,\bolx) \gamma_{ij} \big), 
\end{align}
and  $\Sigma_{ij}$ as in \eqref{defSig}. It  is connected by It\^o's formula to the fully non-linear PDE \eqref{symPDE} via the representation: 
$Y_t$ $=$ $v(t,\boX_t)$, $\boZ_t$ $=$ $D_{\bolx} v(t,\boX_t)$,  $\boG_t^{}$ $=$ $D_{\bolx}^2 v(t,\boX_t)$, $0\leq t\leq T$. 

\vspace{1mm}

Assuming that $G$ is smooth, the algorithm is designed in Algorithm \ref{MKVfullscheme}. 

\begin{small}

\begin{algorithm2e}[!h] \label{AlgofullMKV} 
\DontPrintSemicolon 
\SetAlgoLined 
\vspace{1mm}
{\bf Initialization}: Initialize from the exchangeable function: $\widehat\Uc_{N_T}(\cdot)$ $=$ $G(\cdot)$  and 
the $D$-exchangeable function $\widehat\boZc_{N_T}(\cdot)$ $=$ $D G(\cdot)$. 

\For{ $k$ $=$ $N_T-1,\ldots,0$} 
{minimize over  symmetric NN  $\Uc_k$, and $D$-symmetric NN  $\Zc_k$, the quadratic loss function
\begin{align} \label{loss2}
J_2(\Uc_k,\Zc_k) &= \; \E\Big| \widehat\Uc_{k+1}(\boX_{k+1}^{N_T}) - \Uc_k(\boX_k^{N_T})  \\
&  \quad \quad \quad  + \;   F_{b,\sigma}\big(t_k,\boX_k^{N_T},\Uc_k(\boX_k^{N_T}),\boZc_k(\boX_k^{N_T}),D\widehat\boZc_{k+1}(\boX_{k+1}^{N_T})\big) \Delta t_k  \\
& \quad \quad  \quad - \;  \sum_{i=1}^N \sum_{j=0}^N \big(\Zc_k(\boX_k^{N_T},X_k^{i,N_T})\big) \trans\sigma_{ij}\big(t_k,\boX_k^{N_T}\big) \Delta W_k^j  \Big|^2,
\end{align}
and update $(\widehat\Uc_k,\widehat\Zc_k)$ as the solution to this minimization problem. Here ${\rm D}\widehat\boZc_{k+1}$ is the automatic differentiation of the $D$-symmetric NN  $\widehat\boZc_{k+1}$   
computed previously at the time step $k+1$. 
}
\caption{Fully nonlinear DPBD scheme with symmetric NN  \label{MKVfullscheme} }
\end{algorithm2e}

\end{small}

\subsection{The case of mean-field PDEs}\label{sec: extension}

We consider in this section the case where the PDE \eqref{symPDE} is the particles approximation of a McKean-Vlasov control problem with a Bellman equation \eqref{meanPDE}  in the Wasserstein space of probability measures  as described in Example \ref{exMKV}. 
To simplify the presentation, we  consider that there is only control on the drift coefficient $\beta(t,x,\mu,a)$ but no control on the diffusion coefficient $\sigma(t,x,\mu)$ and $\sigma_0(t,x,\mu)$ of  the McKean-Vlasov equation  \eqref{diffXMKV}.  
In this case, recall that when the solution $v(t,\mu)$ to this Bellman equation is smooth, an optimal control is given in feedback form by (see \cite{phawei17}): 
\beqs
\alpha_t^* &=& 
\hat a(t,X_t^*,\P^0_{_{X_t^*}},\partial_\mu v(t,\P^0_{_{X_t^*}})(X_t^*)),  
\enqs
where $\hat a(t,x,\mu,z)$ is an arg$\min$  of $a$ $\in$ $A$ $\mapsto$ $\beta(t,x,\mu,a).z + f(x,\mu,a)$,  
and $X^*$ $=$ $X^{\alpha^*}$ is the optimal McKean-Vlasov state process.  

\vspace{1mm}

\noindent {\bf Approximation of the optimal control by forward induction of the scheme.} 
As proven in \cite{GPW21b}, the solution $(\boX,Y,\boZ)$  to the FBSDE \eqref{diffX}-\eqref{BSDEsemi} provides an approximation with a rate of convergence $1/N$, when $N$ goes to infinity, of the solution $v$ to \eqref{meanPDE}, and its $L$-derivative: 
$Y_t$ $\simeq$ $v(t,\bar\mu(\boX_t))$, $NZ_t^i$ $\simeq$ $\partial_\mu v(t,\bar\mu(\boX_t))(X_t^i)$.  
The drift coefficients $b_i$ of the forward particles system $\boX$ should be chosen in order to generate from training simulations a  
suitable exploration of the state space for getting a good approximation of the optimal feedback control. 
In practice, in a first step, one can choose $b_i(t,\bolx)$ $=$ $\beta(t,x_i,\bar\mu(\bolx),a_0)$,  for some arbitrary value $a_0$ $\in$ $A$ of the control. 
After a first implementation of Algorithm \ref{MKVscheme}, we thus have an approximation of $\partial_\mu v(t,\mu)(x)$  at time $t$ $=$ $t_k$, and $\mu$ $=$ $\bar\mu(\bolx)$, by $N\widehat\Zc_k(\bolx,x)$. 
Notice however that we solved the PDE along the law of the forward training process, which is different from the optimally controlled process law, except at the initial time $t_0$, where we then get  
an approximation of the optimal feedback control with 
\begin{align}
(x,\bar\mu(\bolx)) & \longmapsto \; \hat a(t_0,x,\bar\mu(\bolx),N\widehat\Zc_0(\bolx,x)).
\end{align}
Next,  by defining an updated initial drift coefficient as
\begin{align}
\hat b_i(t_0,\bolx) & := \; \beta\big(t_0,x_i,\bar\mu(\bolx),\hat a(t_0,x_i,\bar\mu(\bolx),N\widehat\Zc_0(\bolx,x_i) )\big), \; \mbox{ for } \bolx = (x_i)_{i\in \llbracket 1,N\rrbracket}, \; i=1,\ldots,N, 
\end{align}
and considering the $N$-particle discrete-time system $\{\hat\boX_k^{N_T}=(\hat X_k^{i,N_T})_{i \in \llbracket 1,N\rrbracket}, k=0,\ldots,N_T\}$, starting from i.i.d. samples $X_0^i$, $i$ $=$ $1,\ldots,N$ distributed according to some distribution $\mu_0$ on $\R^d$, 
and  with dynamics 
 \begin{align} \label{dynMKVopt}
 \hat X^{i,N_T}_{1} &= \;  X^{i}_{0} +  \hat b_i(t_0,\boX_0) \Delta t_0    +   \sigma(t_0,X_0^i,\bar\mu(\boX_0)) \Delta W_0^i,  \\
\hat X^{i,N_T}_{k+1} &= \;  \hat X^{i,N_T}_{k} +  b_i(t_k,\boX_k^{N_T}) \Delta t_k   +   \sigma(t_k,\hat X^{i,N_T}_{k},\bar\mu(\hat{\boX}_k^{N_T})) \Delta W_k^i, 
 \end{align}
 for $k$ $=$ $1,\ldots,N_T-1$, 
 we obtain an approximation of the distribution of the optimal particle mean-field process at time $t_1$.  
 Applying the algorithm again between $t_1$ and $t_{N_T}$ then allows to compute an approximation of the optimal feedback 
 control $\hat a(t_1,x,\bar\mu(\bolx),\widehat\Zc_1(\bolx,x))$ at time $t_1$ and to update the simulation of $\hat X^{i,N_T}_{2}$. By induction, we can compute the optimal feedback control at every time step through $N_T$ executions of the scheme.
 
 \vspace{1mm}
 
\noindent {\bf Approximation of the solution by randomization of the training simulations.}   
 Algorithm \ref{MKVscheme} provides actually an approximation of $v(t,\mu)$ (resp. $\partial_\mu v(t,\mu)(x)$) at time $t_k$, and for empirical measures $\mu$ $=$ $\bar\mu(\bolx)$, by 
$\widehat\Uc_k(\bolx)$ (resp.  $N\widehat\Zc_k(\bolx,x)$). Thus, in order to get an approximation of $v(t_k,.)$ (resp. $\partial_\mu v(t_k,.)(x)$) on the whole Wasserstein 
space $\Pc_2(\R^d)$, we need a suitable exploration of $\bar\mu(\boX_k^{N_T})$ when using the training simulations $\boX_k^{N_T}$, $k$ $=$ $0,\ldots,N_T$. For that purpose, 
some randomization can first be implemented by randomizing the initial law $\mu_0$ of the forward process. By sampling $\mu_0$ in a compact set $K$ of $\Pc_2(\R^d)$ for each batch element, such as a family of Gaussian measures for instance, our algorithm will be able to learn the value function $v(0,\mu)$ and its Lions derivative $\partial_\mu v(0,\mu)$ on $K$. Therefore, instead of solving the PDE several times for each initial law we can run the algorithm only once. This can be useful if we have an uncertainty in the initial law of the problem we aim to solve. It corresponds to learning the solution $v(t_k,\overline{\mu_{k,\ell}})$ on a family of empirical measures corresponding to forward processes $X^{i,(\ell),N_T}$, $i$ $=$ $1,\ldots,N$,  with initial laws $\mu_0^\ell\in K$. Relying on the generalization properties of neural networks, we expect to approximate the value function at time $t_0$ $=$ $0$ on $K$. 
Furthermore, if the goal is to obtain an approximation of the PDE solution at any time step $t_k$, the task is more complex. 
A randomization needs to be performed at each time step $t_k$ by sampling $\boX^{N_T}_{k}$ according to a Gaussian mixture $\nu_k$ with random parameters. We then apply Algorithm \ref{MKVscheme}, and  expect to learn the solution over measures with regular densities. The updated method is presented in Algorithm \ref{MKVscheme mixture}. If the state space exploration is efficient, the feedback control will be directly available with only one execution of the algorithm, contrarily to the previously described procedure with $N_T$ executions. We should explore the Wasserstein space well enough to learn the value function and its derivative on the unknown law of the optimal process.

\vspace{2mm}

\begin{small}

\begin{algorithm2e}[H]\label{algo mixture}
\DontPrintSemicolon 
\SetAlgoLined 
\vspace{1mm}
\begin{footnotesize}
{\bf Initialization}: Initialize from the exchangeable function: $\widehat\Uc_{N_T}(\cdot)$ $=$ $G(\cdot)$ 

\For{ $k$ $=$ $N_T-1,\ldots,0$} 
{define random variables $$L \sim U(1,L_{max}),\ \varphi_l\sim U(0,1),\ \mu_i\sim U(-\mu_{max},\mu_{max}),\ (\theta_i)^2 \sim U(0,\sigma_{max}^2)$$\\
define a random Gaussian mixture $\nu_k$ of random density $$\frac{\sum_{l=1}^L \varphi_l \Nc(\mu_i,\theta_i^2)}{\sum_{l=1}^L \varphi_l}$$ \\
define $N$ i.i.d. particles $X^{i,N_T}_{k}$ with law $\nu_k$ for $i=1,\cdots,N$, \\
perform one Euler-Maruyama scheme step
\begin{align} \label{mixture}
 X^{i,N_T}_{k+1} &= \;   X^{i,N_T}_{k} +  b_i(t_k,\boX_k^{N_T}) \Delta t_k      +   \sigma\big(t_k, X^{i,N_T}_{k},\bar\mu(\boX_k^{N_T})\big) \Delta W_k^{i}, 
 \end{align}
minimize over  symmetric NN  $\Uc_k$, and $D$-symmetric NN  $\Zc_k$, the quadratic loss function (with $H$ as in \eqref{HMKV}): 
\begin{align} 
J_1(\Uc_k,\Zc_k) &= \; \E\Big| \widehat\Uc_{k+1}(\boX_{k+1}^{N_T}) - \Uc_k(\boX_k^{N_T})  \\
&  \quad \quad \quad  + \;   H\big(t_k,\boX_k^{N_T},\Uc_k(\boX_k^{N_T}),\boZc_k(\boX_k^{N_T})\big) \Delta t_k  \\
& \quad \quad  \quad - \;  \sum_{i=1}^N  \big(\Zc_k(\boX_k^{N_T},X_k^{i,N_T})\big) \trans  \sigma\big(t_k, X^{i,N_T}_{k},\bar\mu(\boX_k^{N_T})\big) \Delta W_k^i  \\
& \quad \quad \quad - \;  \sum_{i=1}^N  \big(\Zc_k(\boX_k^{N_T},X_k^{i,N_T})\big) \trans   \sigma_0\big(t_k, X^{i,N_T}_{k},\bar\mu(\boX_k^{N_T})\big)  \Delta W_k^0   \Big|^2,
\end{align}
and update $(\widehat\Uc_k,\widehat\Zc_k)$ as the solution to this minimization problem. 
}
\end{footnotesize}
\caption{DBDP scheme with symmetric NN and exploration of Wasserstein space \label{MKVscheme mixture} }
\end{algorithm2e}
  
\end{small}

\section{Numerical results} \label{secnum}

In the different test cases,  for the approximation of  the solution $v$  by means of symmetric neutral networks, we used  DeepSets.
 \subsection{A toy example of symmetric PDE in very high dimension}

 We consider a symmetric semi-linear PDE: 
 \begin{align}
 \begin{cases}
 \partial_t v + b.D_\bolx v + \frac{1}{2}{\rm tr}(\sigma\sigma\trans D_\bolx^2 v)  + f(\bolx,v,\sigma\trans D_\bolx v) \; = \; 0, \quad (t,\bolx) \in [0,T)\times\R^N \\
 v(T,\bolx) \; = \; \cos(\bar x), \quad  \mbox{ with } \bar x = \sum_{i=1}^N x_i, \mbox{ for }   \bolx = (x_1,\ldots,x_N) \in \R^N, 
 \end{cases}
 \end{align}
 with $b$ $=$ $0.2/N$, $\sigma$ $=$ $\frac{I_N}{\sqrt{N}}$, 
 \begin{align}
 f(\bolx,y,\boz) = \left(\cos(\overline{x})+0.2\sin(\overline{x})\right)e^{\frac{T-t}{2}} - \frac{1}{2} (\sin(\overline{x})\cos(\overline{x})e^{T-t})^2+ \frac{1}{2N}(y(1_N \cdot \boz))^2.
\end{align}
so  that the PDE solution is exchangeable and given by
\begin{align}
v(t,\bolx) &= \;  \cos\left(\overline{x}\right) \exp\big(\frac{T-t}{2}\big).
\end{align}
We solve this PDE in  dimension $N$ $=$ $1000$ by using the deep backward scheme  (DBDP)  in \cite{hure2020deep} with  $60$ time steps, and estimate $U_0$ $=$ $v(0,1_N)$ and $Z_0$ $=$ $D_\bolx v(0, 1_N)$.  
For the approximation of $v$, and its gradient $D_\bolx v$, we test with three classes of networks: 
\begin{itemize}
\item[(i)] DeepSet $\Uc$ for $v$, and AD-DeepSet ${\rm D}\Uc$  for $D_\bolx v$ (DeepSets derivative case).  
\item[(ii)] DeepSet  for $v$, and DeepDerSet for $D_\bolx v$ (DeepDerSet case)  
\item[(iii)] Feedforward  for $v$ and $D_\bolx v$ (Feedforward case) 
\end{itemize}
For each of theses case, we use ReLU activation functions for all the networks, and  for the feedforward network, we choose 3 layers of $1010$ neurons.

\begin{Remark}
{\rm 
An alternative to Case (i) is to consider an  AD-DeepSet  ${\rm D}\Uc_1$ for $D_\bolx v$ with  $\Uc_1$ another DeepSet independent of the one $\Uc$ used for $v$.
}
\ep
\end{Remark}

We report the solution in Table \ref{tab:pdeCost}.

\begin{table}[h!]
    \centering
    \begin{tabular}{|c|c|c|c|c|c|c|c|} \hline
       \multicolumn{2}{|c|}{ Analytical }   & \multicolumn{2}{|c|}{(i) AD-DeepSets }   &
       \multicolumn{2}{|c|}{(ii) DeepDerSet} & \multicolumn{2}{|c|}{(iii) Feedforward}  \\ \hline
        $U_0$ & $Z_0$ & $U_0$ & $Z_0$ & $U_0$ & $Z_0$ &$U_0$ & $Z_0$ \\ \hline
            0.9272  &  -1.3632     &  0.9289   & -1.2973  &0.90140   & -1.304  & 0.6896 &  -1e-7 \\ \hline
    \end{tabular}
    \caption{PDE resolution in dimension 1000 with DBDP scheme \cite{hure2020deep}. }
    \label{tab:pdeCost}
\end{table}

We observe that  the results with the feedforward network are  not good. This is due to the fact that the feedforward network is not able to approximate 
correctly the final condition whatever the initial learning rate and the number of epochExt are taken, as already shown in Table \ref{tab:func3}.  In contrast, we see that 
the AD-DeepSets and DeepDerSet networks give  good results but only when the initial learning rate is taken small enough (here we took  \num{1e-5}). 
Finally, we  have tested the Deep BSDE method in \cite{han2018solving} with the variation proposed in \cite{CWNMW19}  using a network reported in section \ref{sec:tAndX}. The results are unstable and so we do not report them.  
A direct use of \cite{han2018solving} method with a network per time step is impossible to test due the size of the problem but results in lower dimension also indicate some instability directly linked to the initialization of the network.

\subsection{A mean-field control problem of systemic risk} \label{subsecex1}

We consider a mean-field model of systemic risk introduced in \cite{carfousun}.  This model was introduced in the context of mean field games but here we consider a cooperative version. The limiting pro\-blem (when the number of banks is large) of the social planner (central bank) is formulated as follows. The log-monetary reserve of the representative bank is governed by the mean-reverting controlled McKean-Vlasov dynamics 
\begin{align}
dX_t &= \; \big[  \kappa( \E[X_t] - X_t) + \alpha_t]  \ \di t +  \sigma d W_t , 
\quad X_0 \sim \mu_0, 
\end{align}
where $\alpha$ $=$ $(\alpha_t)_t$  is the  control rate of borrowing/lending to a central bank that aims to minimize  the functional cost
\begin{align} \label{defV0}
J(\alpha) \; = \; \E \Big[ \int_0^T \tilde f(X_t,\E[X_t],\alpha_t) \ \di t + \tilde g(X_T,\E[X_T]) \Big] &\;  \rightarrow \quad V_0 \; = \; \inf_{\alpha} J(\alpha), 
\end{align}
where the running and terminal costs are given by 
\begin{align}
\tilde f(x,\bar x,a)  \; = \; \frac{1}{2}a^2 - qa(\bar x-  x) + \frac{\eta}{2} (\bar x-x)^2, & \quad \tilde g(x,\bar x) \; = \; \frac{c}{2}(x-\bar x)^2, 
\end{align} 
for some positive constants $q$, $\eta$, $c$ $>$ $0$, with $q^2$ $\leq$ $\eta$.  

The value function $v$  to the mean-field type control problem \eqref{defV0} is solution to the Bellman (semi-linear PDE) equation \eqref{meanPDE} with $\sigma$ constant, $\sigma_0$ $\equiv$ $0$, $r$ $=$ $0$, and
\begin{align}
h(t,x,\mu,z,\gamma) 
& = \; \inf_{a \in \R} \big\{  \big[\kappa( \E_\mu[\xi] - x) + a\big]z   + \frac{1}{2}a^2 -qa(\E_\mu[\xi] -x) \big\} + \frac{\sigma^2}{2}  \gamma + \frac{\eta}{2} \big(\E_\mu[\xi] -x\big)^2 \\
& = \; (\kappa + q)( \E_\mu[\xi] - x)z +  \frac{\sigma^2}{2}  \gamma   +  \frac{\eta-q^2}{2} \big(\E_\mu[\xi]-x\big)^2 - \frac{z^2}{2},
\end{align}
and $\Gc(\mu)$ $=$ $\frac{c}{2}{\rm Var}(\mu)$ $:=$ $\frac{c}{2}\E_\mu | \xi - \E_\mu[\xi]|^2$ is  the variance of the distribution $\mu$ (up to $c/2$). Here, we use the notation: $\E_\mu[\varphi(\xi)]$ $:=$ $\int \varphi(x) \mu(dx)$. 

The finite-dimensional approximation of \eqref{defV0} with $N$-bank model corresponds to the symmetric  Bellman semi-linear PDE on $[0,T]\times\R^N$: 
\begin{align}\label{eq: HJB SYSTEMIC}
\partial_t v^N + \sum_{i=1}^N (\kappa + q)\big(\bar x  - x_i\big)\partial_{x_i} v^N + \frac{\sigma^2}{2}  \Delta_{\boldsymbol{x}} v^N 
+  \frac{\eta - q^2}{2N}\sum_{i=1}^N \big(\bar x  - x_i \big)^2  -  \frac{N}{2}\sum_{i=1}^N |\partial_{x_i} v^N|^2 & = 0,
\end{align}
for $\boldsymbol{x}=(x_1,\ldots,x_N)$ $\in$ $\R^N$, where we set $\bar x$ $=$ $\frac{1}{N}\sum_{i=1}^N x_i$, and $\Delta_{\boldsymbol{x}}$ $=$ $\sum_{i=1}^N \partial_{x_ix_i}^2$ is the Laplacian operator. 
We numerically solve \eqref{eq: HJB SYSTEMIC} with Algorithm \ref{MKVscheme}  described in Section \ref{secalgo}. The algorithm is trained with the forward 
 process in $\R^N$: 
\begin{align}
X_{k+1}^i & = \; X_k^i  + \sigma \Delta W_k^i , \quad X_0^i \sim \mu_0, \;\; k=0,\ldots,N_T -1, \;\;\; i=1,\ldots,N.  
\end{align}     
The choice of a null  drift for this training process is intuitively justified by the fact that the objective  in \eqref{defV0} is to  incite the log-monetary reserve of the banks to be close to the average of all the other banks, hence we formally expect their drift 
to be close to zero.  

We test our algorithm by increasing $N$, and compare with the explicit solution of the limiting linear-quadratic McKean-Vlasov control problem \eqref{defV0}, which is solved via the resolution of a Riccati equation (see \cite{baspha19}), and is analytically given by 
\begin{align}
v(t,\mu) &= \; K_t \mathrm{Var}(\mu) + \sigma^2 \int_t^T K_s \ \di s,
\end{align}
where 
\begin{align}
K_t & = 
- \frac{1}{2} \Big[ \kappa + q - \sqrt{\Delta} \frac{ \sqrt{\Delta} \sinh(\sqrt{\Delta}(T-t))  + (\kappa + q + c) \cosh (\sqrt{\Delta}(T-t))}{ \sqrt{\Delta} \cosh(\sqrt{\Delta}(T-t))  + 
(\kappa + q + c) \sinh (\sqrt{\Delta}(T-t))} \Big],   
\end{align}     
with $\sqrt{\Delta}$ $=$ $\sqrt{(\kappa+q)^2  + \eta -  q^2}$,  and 
\begin{align}
\int_t^T K_s \ \di s & = \;  \frac{1}{2} \ln \Big[ \cosh (\sqrt{\Delta}(T-t)) + \frac{\kappa + q  + c}{\sqrt{\Delta}}  \sinh (\sqrt{\Delta}(T-t)) \Big] - \frac{1}{2} (\kappa + q)(T-t). 
\end{align}

We have tested various approximation symmetric networks and different resolution me\-thods. We list the methods that fail to solve the problem   in high dimension:   
\begin{itemize}
    \item First we tried the global resolution method in \cite{weinan2017deep} by using the network described in paragraph  \ref{sec:tAndX} for $v$. 
    
    In this case, we could not  obtain exploitable results. 
    \item Then we decided to use the local resolution method \cite{hure2020deep} with a DeepDerSet approximation approach  for $\Zc$. We found that the method give accurate results  in dimension below 100  but with a variance increasing with the dimension. Results were impossible to exploit in dimension 1000. We thus  decided to not report the results.
    \item At last we tested the local resolution methods  \cite{hure2020deep}  with  classical feedforward networks using tanh or ReLU activation functions. Two variants were tested with $N=500$: the first one using a network for $v^N$ (we stress the dependence of the solution to the PDE on $N$) and another network for $Dv^N$ giving values not exploitable, and a second version using a single network for $v^N$ and using automatic differentiation to approximate $Dv^N$ giving a very high bias and a high standard deviation.
\end{itemize}
Therefore we only report the case where we use a  DeepSet network for $\Uc$ and a second AD-DeepSet  network  to estimate  $\Zc$ or  a  single DeepSet network for $\Uc$ which is differentiated to approximate $\Zc$.

We test the  tanh and ReLU activation function on this test case using the parameters
$\sigma=1$, $\kappa=0.6$, $q=0.8$, $c=2$, $\eta=2$, $T=1$. We report  $v^N$ estimated with different values of $N_T$ and $N$ at time $t$ $=$ $0$ and $x=0$ so using $\mu_0=\delta_0$ on the figures below. The theoretical solution obtained  when $N$ goes to infinity is $0.29244$.\\
We use a batch size equal to $200$, a number of gradient iteration equal to $30000$ for the resolution to project the terminal condition on the network and $6000$  gradient iterations for other resolutions. The initial learning rate is taken equal to $1e-4$ at the first resolution and $5e-5$ for other resolutions. The learning is taken decreasing linearly with gradient iterations to $5e-6$.

On figure \ref{fig:systemicRelUUZ}, we give the results obtained with  ReLU activation function using
a  DeepSet network for $\Uc$ and a second AD-DeepSet network to approximate $\Zc$. The convergence is steady as $N_T$ grows and as the dimension grows leading to a very accurate result for $N_T=60$ and $N=1000$.
\begin{table}[H]
    \centering
    \begin{tabular}{|c|c|c|c|c|}
        \hline
 $N_T$ &      Dimension $N$ & Averaged  & Std & Relative error\\
        \hline
  15   &   10 & 0.259&  0.0029& 0.11 \\ \hline  
  15   &   100 & 0.2871  & 0.0016 & 0.018 \\ \hline  
  15   &   500 & 0.2866& 0.00179 &  0.019\\ \hline 
  15   &   1000 & 0.2877 &  0.00238 &  0.016\\ \hline 
  30   &   10 &  0.265 & 0.004  &  0.09  \\ \hline  
  30   &   100 & 0.2892  & 0.2892 &  0.010\\ \hline  
  30   &   500 &  0.2897 & 0.00153& 0.009 \\ \hline 
  30   &   1000 & 0.2899&  0.00146 & 0.0084\\ \hline 
  60   &   10 & 0.2655 & 0.0045  & 0.092 \\ \hline  
  60   &   100 & 0.2894& 0.0012 &  0.010\\ \hline  
  60   &   500 & 0.2894 & 0.0027& 0.010 \\ \hline 
  60   &   1000 & 0.2916 & 0.0014& 0.0025 \\ \hline   
    \end{tabular}
    \caption{Systemic risk with ReLU activation function, a  DeepSet network for $\Uc$ and a second AD-DeepSet  network  to estimate  $\Zc$.}
    \label{fig:systemicRelUUZ}
\end{table}

Using a ReLU activation function, a single network for $\Uc$ which is differentiated to approximate $\Zc$, we get the results in figure \ref{fig:systemicRelUUDU}. The convergence is still steady but results are not as good as in table \ref{fig:systemicRelUUZ}.

\begin{figure}[H]
    \centering
    \begin{minipage}[c]{.32\linewidth}
    \includegraphics[width=\linewidth]{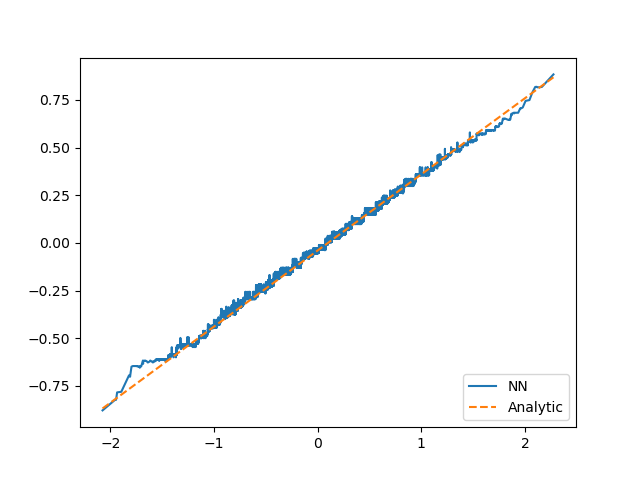}
    \caption*{$N_T=10$}
    \end{minipage}
    \begin{minipage}[c]{.32\linewidth}
    \includegraphics[width=\linewidth]{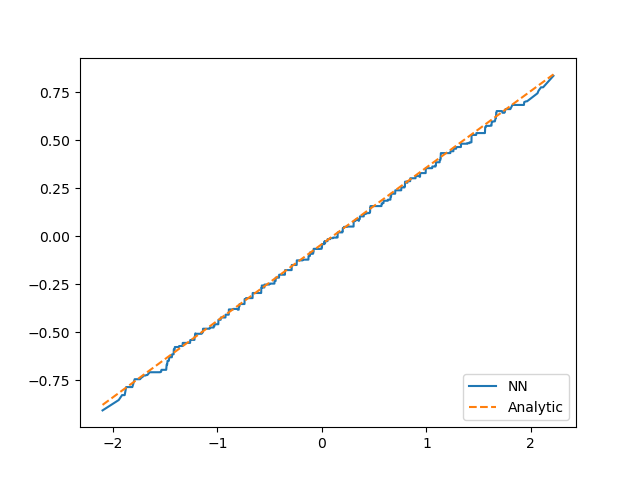}
    \caption*{$N_T=20$}
    \end{minipage}
    \begin{minipage}[c]{.32\linewidth}
    \includegraphics[width=\linewidth]{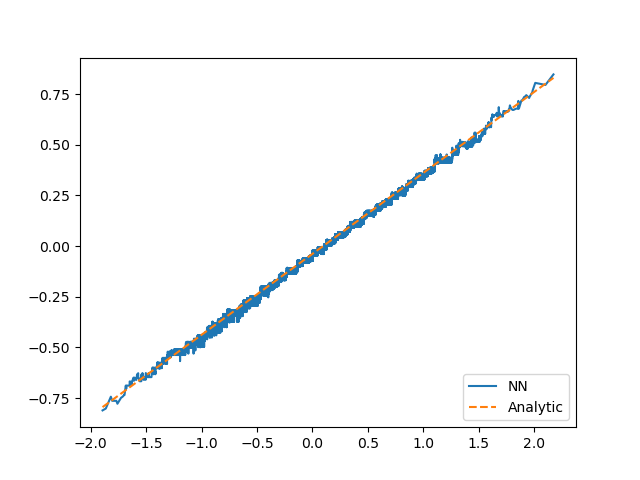}
    \caption*{$N_T=40$}
    \end{minipage}
    \caption{Resolution on $[0.5,1]$  in dimension $N=500$ : analytic Lions derivatives versus $N\Zc$ estimated by the network. DeepSet network for $\Uc$, AD-DeepSet for $\Zc$. ReLU activation function.}
    \label{fig:LionSDeriv500}
\end{figure}

\begin{table}[H]
    \centering
    \begin{tabular}{|c|c|c|c|c|}
        \hline
 $N_T$ &      Dimension $N$ & Averaged  & Std & Relative error \\
        \hline
  15   &   10 & 0.2530& 0.0074 &  0.1346 \\ \hline  
  15   &   100 & 0.27968  & 0.0051 &  0.043\\ \hline  
  15   &   500 &0.2938 & 0.0067&  0.0049\\ \hline 
  15   &   1000 & 0.3084& 0.0253 & 0.054 \\ \hline 
  30   &   10 & 0.2494 & 0.0074 & 0.1471   \\ \hline  
  30   &   100 & 0.2756&  0.00677& 0.057 \\ \hline  
  30   &   500 & 0.2885& 0.0127& 0.013 \\ \hline 
  30   &   1000 & 0.2860& 0.009 & 0.02 \\ \hline 
  60   &   10 & 0.2519  & 0.0037 &  0.138 \\ \hline  
  60   &   100 & 0.28253 & 0.0047 & 0.033 \\ \hline  
  60   &   500 & 0.28329& 0.0108& 0.03 \\ \hline 
  60   &   1000 &0.2881 & 0.0043&  0.014\\ \hline 
    \end{tabular}
    \caption{Systemic risk  with ReLU activation function, a  single DeepSet network for $\Uc$ which is differentiated to approximate $\Zc$. }
    \label{fig:systemicRelUUDU}
\end{table}
The replace the ReLU activation function by a tanh one  using two networks and the results are given in table~\ref{fig:systemicTanhUZ}. The convergence is not steady and increasing to much $N$ or $N_T$ worsen to results : it shows the importance of the activation function in this method.
\begin{table}[H]
    \centering
    \begin{tabular}{|c|c|c|c|c|}
        \hline
 $N_T$ &      $N$ & Averaged  & Std & Relative error \\
        \hline
  15   &   10 &  0.2678   & 0.0061  & 0.08\\ \hline  
  15   &   100 & 0.28858  &  0.0144 &  0.013\\ \hline  
  15   &   500 & 0.2491 & 0.027 &  0.14\\ \hline 
  15   &    1000& 0.27401 &  0.0127 &  0.063\\ \hline 
  30   &    10&    0.2725 & 0.0052 & 0.068  \\ \hline  
  30   &   100 & 0.2959  &  0.0161 & 0.012 \\ \hline  
  30   &   500 & 0.2577  & 0.01568&  0.118\\ \hline 
  30   &   1000 & 0.320 &  0.0030 & 0.096\\ \hline 
  60   &   10 & 0.2739& 0.0049  & 0.063 \\ \hline  
  60   &   100 & 0.2924&  0.0309& 0.0001\\ \hline  
  60   &   500 & 0.3158& 0.00297& 0.079 \\ \hline 
  60   &   1000 & 0.2210 &0.004 & 0.24 \\ \hline 
    \end{tabular}
    \caption{ Systemic  with tanh activation function, a  DeepSet network for $\Uc$ and a second AD-DeepSet  network  to estimate  $\Zc$.}
    \label{fig:systemicTanhUZ}
\end{table}
At last we do not report the test obtained using a ReLU activation function  for the first network and a tanh one for the second network given results far better than in table \ref{fig:systemicTanhUZ} but not as good as in tables \ref{fig:systemicRelUUZ} and \ref{fig:systemicRelUUDU}. 
We also test the accuracy of our algorithm for approximating the $L$-derivative of the solution, which is here explicitly given by 
\begin{align}
\partial_\mu v(t,\mu)(x) &= 2  K_t (x - \E_\mu[\xi]). 
\end{align}
For this purpose, using $N_T$ steps,  we solve the same problem  on $[t,T]$,  starting at $t=\frac{T}{2}$ with a distribution $\mu_0$ equal to real distribution of the solution of \eqref{defV0} taken at date $t$.  After training, we plot $x$ $\mapsto$ $N\widehat\Zc(\boX_t,x)$, where $\boX_t$ $\sim$ $\mu_0^{\otimes N}$, and compare to the analytic solution: 
$x$ $\mapsto$ $\partial_\mu v(t,\mu_0)(x)$. 
Some graphs are reported in figure  \ref{fig:LionSDeriv500}, which shows the accuracy of the approximation. 

\begin{figure}[H]
    \centering   
    \begin{minipage}[c]{.49\linewidth}
    \includegraphics[width= 0.8\linewidth]{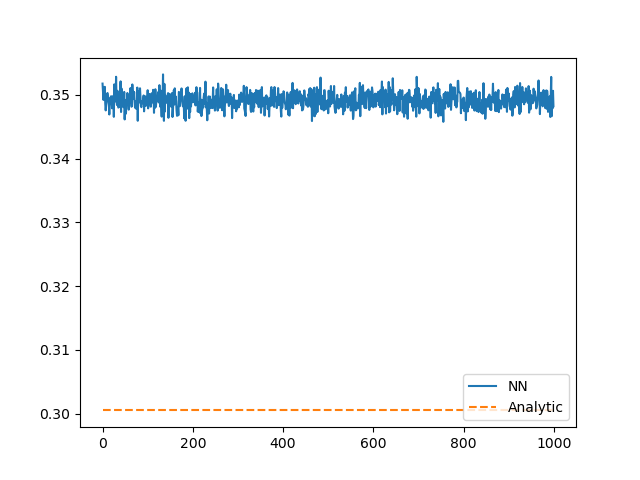}
    \caption*{Solution for $\mu_0 = \Nc(-0.8,0.09)$}
    \end{minipage}
    \begin{minipage}[c]{.49\linewidth}
    \includegraphics[width= 0.8\linewidth]{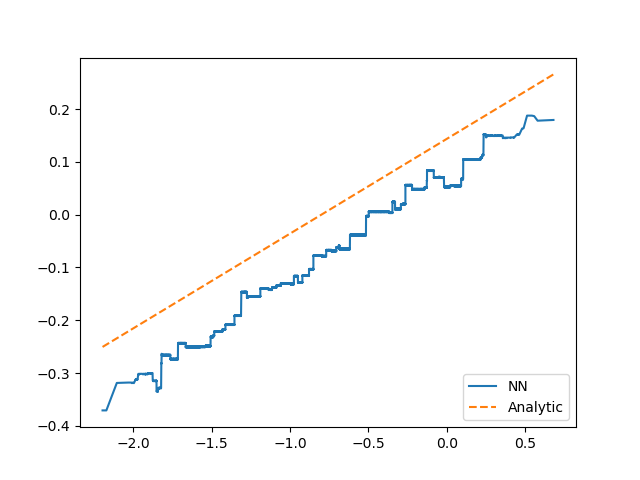}
    \caption*{Lions derivative for $\mu_0 = \Nc(-0.8,0.09)$}
    \end{minipage}
        \begin{minipage}[c]{.49\linewidth}
    \includegraphics[width= 0.8\linewidth]{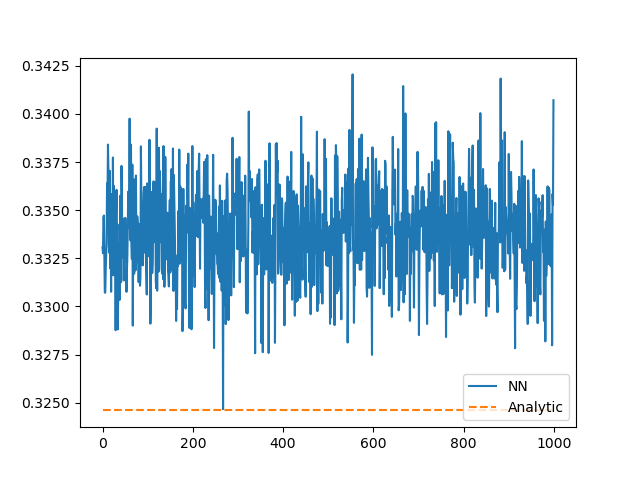}
    \caption*{Solution for $\mu_0 = \Nc(0.,0.36)$}
    \end{minipage}
    \begin{minipage}[c]{.49\linewidth}
    \includegraphics[width= 0.8\linewidth]{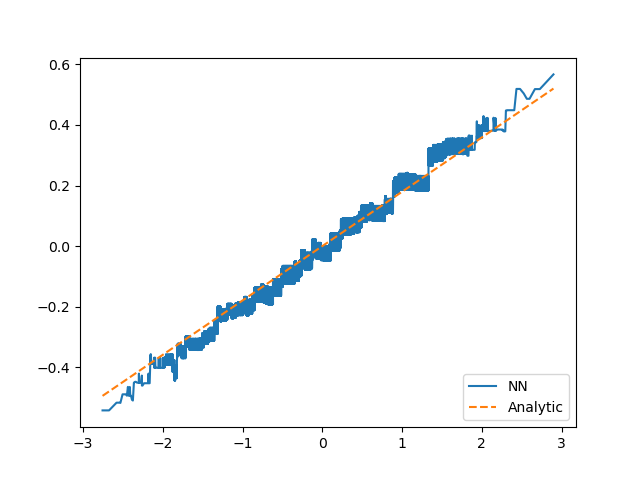}
    \caption*{Lions derivative for $\mu_0 = \Nc(0.,0.36)$}
    \end{minipage}
    \begin{minipage}[c]{.49\linewidth}
    \includegraphics[width= 0.8\linewidth]{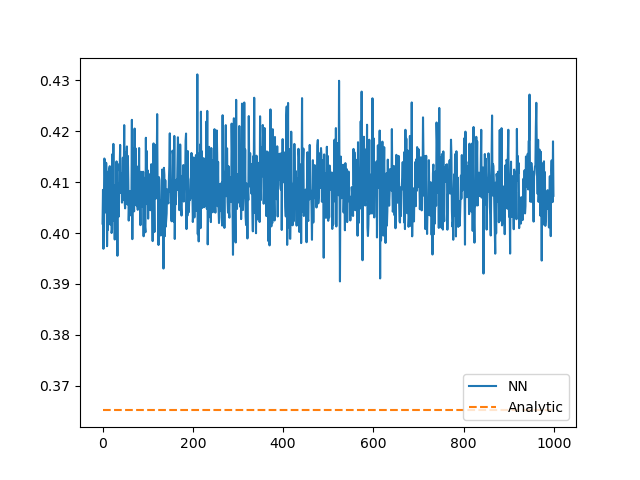}
    \caption*{Solution for $\mu_0 = \Nc(0.8,0.81)$}
    \end{minipage}
    \begin{minipage}[c]{.49\linewidth}
    \includegraphics[width= 0.8\linewidth]{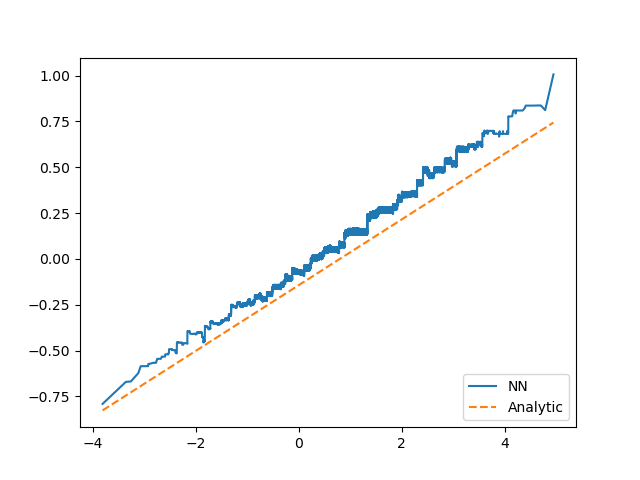}
    \caption*{Lions derivative for $\mu_0 = \Nc(0.8,0.81)$}
    \end{minipage}

    \caption{Solution and Lions derivative after a single training,  with $N=500$, $N_T=30$,  with ReLU activation function, a  single DeepSet network for $\Uc$ which is differentiated to approximate $\Zc$. For the solution, the $x$-axis  corresponds to the sample number and the $y$-axis is the value of the estimated solution. For the Lions derivative, the $x$-axis is the state space and the $y$-axis is the value of the derivative.}
    \label{fig:sysSingelTraining}
\end{figure}

As mentioned in Section \ref{sec: extension}, in theory, the proposed methodology should learn the solution for any initial law $\mu_0$ in the space of measures so that we should be able to solve the problem in infinite dimension. We test our algorithm by sampling $\mu_0$ in the following way: for a sample $j$, we pick up a mean $\hat M \in [-1,1]$ and a standard deviation $\sigma \in [0.2,1]$ with an uniform law.
Then $X^{i,j}_0 \sim \Nc(\hat M,  \sigma^2) \;\;\; i=1,\ldots,N $ and as before we use the forward process:
\begin{align}
X_{k+1}^{i,j} & = \; X_k^{i,j}  + \sigma \Delta W_k^{i,j} , \;\; k=0,\ldots,N_T -1, \;\;\; i=1,\ldots,N.  
\end{align}   
After the training part,  we try to recover the initial solution and the initial Lions derivative for a given $\mu_0$ following a gaussian law.
Results are given on figure \ref{fig:sysSingelTraining}. The Lions derivative is relatively correctly calculated but the initial value can get an error around $15 \%$.

\begin{figure}[H]
    \centering
     \begin{minipage}[c]{.49\linewidth}
    \includegraphics[width=\linewidth]{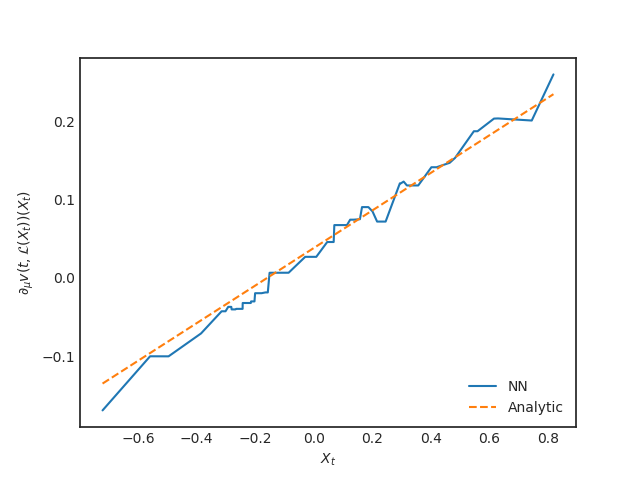}
    \caption*{$t=0.2$}
    \end{minipage}
     \begin{minipage}[c]{.49\linewidth}
    \includegraphics[width=\linewidth]{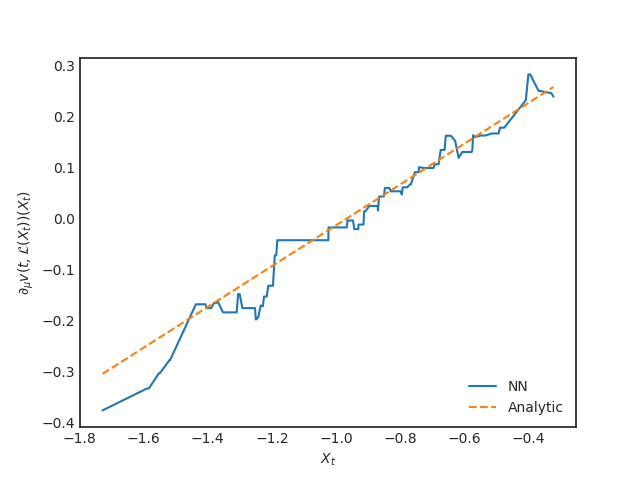}
    \caption*{$t=0.5$}
    \end{minipage}
     \begin{minipage}[c]{.49\linewidth}
    \includegraphics[width=\linewidth]{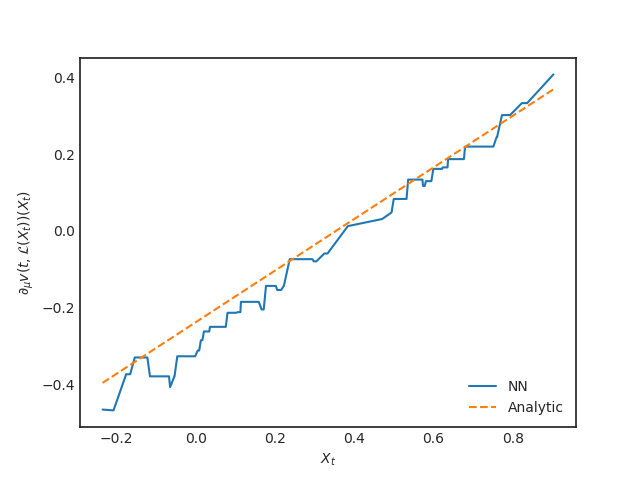}
    \caption*{$t=0.7$}
    \end{minipage}
   \caption{
    Analytic Lions derivative versus $N\Zc$ estimated by the network. Dimension $N=300$, number of time steps $N_T=30$. We use a DeepSet network for $\Uc$ with ReLU activation functions, and $\Zc$ its automatic derivative.
    }
    \label{fig:Lderiv}
\end{figure}

More generally, if we want to solve the PDE at each time step in the Wasserstein space, we can use Algorithm \ref{algo mixture}. In order to illustrate the exploration of the Wasserstein space we  plot in Figure \ref{fig:Lderiv}  the graphs of  $(X^i, N \hat\Zc(\boX,X^i))$, $i$ $=$ $1,\ldots,N,$ vs $X^i$ $\mapsto$ $\partial_\mu v(t,\P_{X^i})(X^i)$,  when $X_t^i$ $\leadsto$ random mixture of Gaussian laws, for $N$ $=$ $300$, $N_T$ $=$ $30$. We observe that we are able to estimate correctly the Lions derivative of the solution (and therefore the optimal control) on several probability measures through a randomized training. Concerning the solution itself, we observe similar behavior as in Figure \ref{fig:sysSingelTraining} with an error of order 10-15\% so we do not show the plots. Further numerical studies are left to future research to improve the estimation of the solution with the randomization procedure. 

\subsection{Mean-variance problem}
\label{sec:MeanVar}

We consider the celebrated Markowitz portfolio selection problem where an investor can invest at any time $t$ an amount $\alpha_t$ in a risky asset (assumed for simplicity to follow a Black-Scholes model with constant rate of return $\beta$ and volatility $\nu$ $>$ $0$), hence generating a wealth process $X$ $=$ $X^\alpha$ with dynamics
\begin{align*}
dX_t &= \; \alpha_t \beta dt + \alpha_t \nu dW_t, \quad 0 \leq t \leq T, \; X_0 = x_0 \in \R.  
\end{align*}
The goal is then to minimize over portfolio control $\alpha$ the mean-variance criterion: 
\begin{align*}
J(\alpha) &= \; \lambda {\rm Var}(X_T^\alpha) - \E[X_T^\alpha],
\end{align*}
 where $\lambda$ $>$ $0$ is a parameter related to the risk aversion of the investor.  Due to the presence of the variance term ${\rm Var}$ in the criterion, the Markowitz problem falls into the class of McKean-Vlasov control problems, and  the associated value function $v$ satisfies the  Bellman equation \eqref{meanPDE} on $[0,T]\times\Pc_2(\R)$ with $r$ $=$ $0$, $\sigma_0$ $\equiv$ $0$, 
\begin{equation}
\begin{cases}
h(x,\mu,z,\gamma) \; = \; \inf_{a\in\R}\big[ z a\beta + \frac{1}{2} \gamma a^2 \nu^2 \big] \; = \;   - \frac{R}{2} \frac{z^2}{\gamma}, \quad z \in \R, \gamma > 0,    \\
\Gc(\mu) \; = \;  \lambda  \E_\mu \big|\xi - \E_\mu[\xi]\big|^2 - \E_\mu[\xi], \quad \mu \in \Pc_2(\R), 
\end{cases}
\end{equation} 
where we set $R$ $:=$ $\beta^2/\nu^2$. 

The associated finite-dimensional PDE with $N$ particles is given by 
\begin{equation} \label{PDE mean variance} 
\begin{cases}
 \partial_t v^N 
  - \frac{R}{2}  \Sum_{i=1}^N \frac{(D_{x_i} v^N)^2}{D_{x_i}^2 v^N}  =   0, \quad  t \in [0,T),  \bolx = (x_1,\ldots,x_N) \in  (\R^d)^N, \\ 
  v^N(T,\bolx) \; = \; \Gc(\bar\mu(\bolx)).
\end{cases}
\end{equation}
We refer to \cite{ismpha19} for the McKean-Vlasov approach to Markowitz mean-variance problems (in a more general context), and we recall that the solution to the Bellman equation is given by  
\begin{align}
\label{eq:analMV}
& v(t,\mu)  = \;  \lambda e^{-R(T-t)}  \E_\mu \big|\xi - \E_\mu[\xi]\big|^2 - \E_\mu[\xi]  - \frac{1}{4\lambda} \big[ e^{R(T-t)} - 1\big] \\
& \partial_\mu v(t,\mu)(x)  = \; 2 \lambda e^{-R(T-t)}  (x - \E_\mu[\xi]) - 1, \quad  \partial_x \partial_\mu v(t,\mu)(x)  = \; 2\lambda e^{-R(T-t)}
\end{align}
and in particular $V_0$ $:=$ $\inf_\alpha J(\alpha)$ $=$ $v(0,\delta_{x_0})$ $=$ $-x_0 - \frac{1}{4\lambda}[e^{RT} - 1]$.  Moreover,  the optimal portfolio strategy is given by 
\begin{align}
\alpha_t^* \; = \; \hat a(t,X_t^*,\E[X_t^*]) & := \; - \frac{\beta}{\nu^2} \Big[ X_t^* -  \E[X_t^*] - \frac{e^{R(T-t)}}{2\lambda} \Big] \\
& = \; - \frac{\beta}{\nu^2} \Big[ X_t^* -  x_0  - \frac{e^{RT}}{2\lambda} \Big], \quad 0 \leq t \leq T, 
\label{eq:optCont}
\end{align}
where $X^*$ $=$ $X^{\alpha^*}$ is the optimal wealth process. 

We test our Algorithm \ref{MKVfullscheme} described in Section \ref{secfull}  with the training of the forward process 
\begin{align}
X_{k+1}^{i,N,\pi} & = \; X_k^{i,N,\pi}  + \frac{R}{2\lambda} \Delta t_k  +  \frac{\sqrt{R}}{2\lambda}  \Delta W_k^i, \quad X_0^i = x_0, \;\; k=0,\ldots,N_T -1, \;\;\; i=1,\ldots,N,  
\end{align}   
which is the time discretization of the  wealth process for a constant portfolio strategy $\alpha_t$ $=$ $\beta/(2\nu^2\lambda)$, which is known to be optimal for the exponential utility function 
$U(x)$ $=$ $-e^{-2\lambda x}$.  This corresponds to the choice of $b_i$ $=$ $R/(2\lambda)$ and $\sigma_{ij}$ $=$ $\sqrt{R}/(2\lambda)$.  Here, notice that 
$\partial_\mu G(\mu)(x)$ $=$ $2\lambda (x - \E_\mu[\xi]) - 1$,  and we then use for the initialization at terminal step $N_T$, the  DeepDerSet  function 
$\Zc_{N_T}((x_i)_i,x)$ $=$ $2\lambda\big(x - \frac{1}{N}\sum_i x_i \big) - 1$ (corresponding to the average function $\mfs((x_i)_i)$ $=$ $\frac{1}{N}\sum_i x_i$),  which yields the automatic differentiation 
${\rm D}\Zc_{N_T}((x_i)_i,x)$ $=$ $2\lambda(1-\frac{1}{N})$. 


\vspace{1mm}

We choose the parameters $\beta = 0.15$, $ \nu = 0.35$, $\lambda = 1$, and the quantile at $99.9\%$ for the truncation in scheme \cite{phawarger20}, and report the results in Table \ref{tab:MV}. The optimization parameters are the same as in the semi linear case, except the batch size taken equal to $50$ and the number of gradient iterations after first step taken equal to $4000$. We use a ReLU Deepset for $\Uc$ and a AD-Deepset  with a tanh activation function for $\Zc$. Remark that in this case it is not possible to use a ReLU activation function for the second network.

\begin{table}[H]
    \centering
    \begin{tabular}{|c|c|c|c|c|}
        \hline
 $N_T$ &      Dimension $N$ & Averaged  & Std & Relative error \\
        \hline
  10   &   10 & -1.0561 & 0.001 &  0.005 \\ \hline  
  10   &   100 &  -1.0522 & 0.0008 & 0.0017 \\ \hline  
  20   &   10 & -1.0570 & 0.0008 &  0.006 \\ \hline  
  20  &     100 & -1.0520& 0.0007& 0.0015 \\ \hline 
  30   &   10 & -1.0578 & 0.0011 &  0.007 \\ \hline  
  30   &   100 & -1.0535& 0.0021 & 0.0029 \\ \hline  
   \end{tabular}
    \caption{\small{Estimate of $\E\big[ v^N(0,X_0^1,\ldots,X_0^N) \big]$ with a deterministic initial condition $X_0 = 1,\ T=1,\ \sigma = 1$. Average and standard deviation
observed over 10 independent runs are reported. The theoretical solution is $-1.0504058$ when $N,N_T \rightarrow + \infty$.}}
\label{tab:MV}
\end{table}

\begin{figure}[H]
    \centering   \begin{minipage}[c]{.38\linewidth}
    \includegraphics[width=0.90\linewidth]{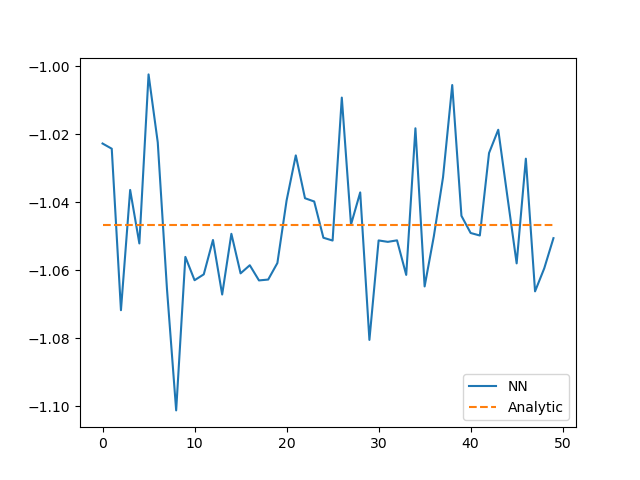}
    \caption*{Solution}
    \end{minipage}
        \begin{minipage}[c]{.38\linewidth}
    \includegraphics[width=0.90\linewidth]{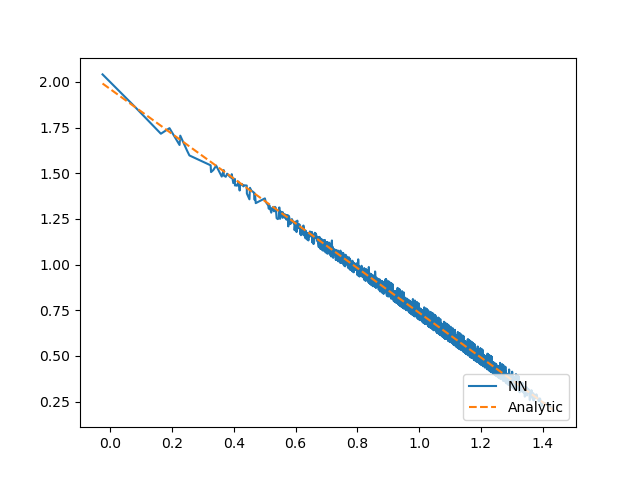}
    \caption*{Control}
    \end{minipage}
    \begin{minipage}[c]{.38\linewidth}
    \includegraphics[width=0.90\linewidth]{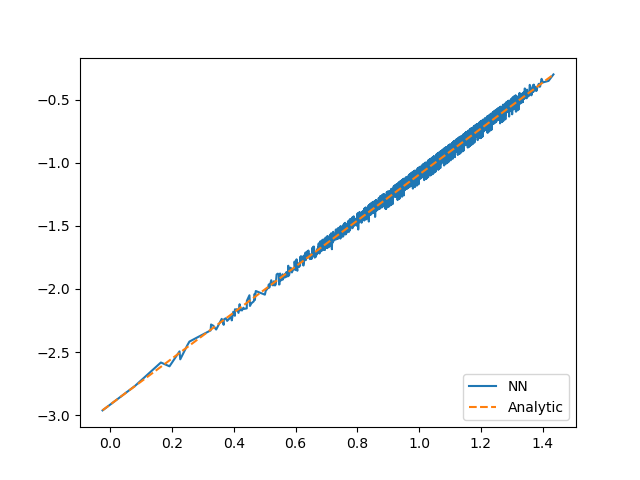}
    \caption*{Derivative}
    \end{minipage}
     \begin{minipage}[c]{.38\linewidth}
    \includegraphics[width=0.90\linewidth]{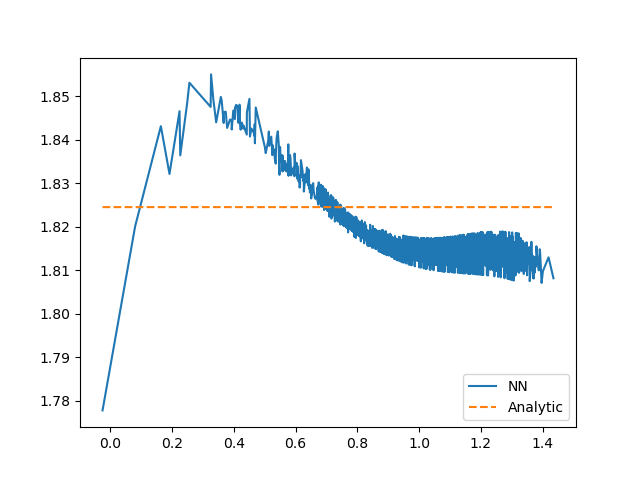}
    \caption*{Second order derivative}
    \end{minipage}
    \caption{Solution and control obtained on the mean variance case at $\frac{T}{2}$ in dimension $100$ with $20$ time steps comparing analytic solution to the calculated one (NN). Truncation factor equal to $0.999$.}
    \label{fig:solAnalMV100}
\end{figure}

\begin{figure}[H]
    \centering   \begin{minipage}[c]{.38\linewidth}
    \includegraphics[width=0.90\linewidth]{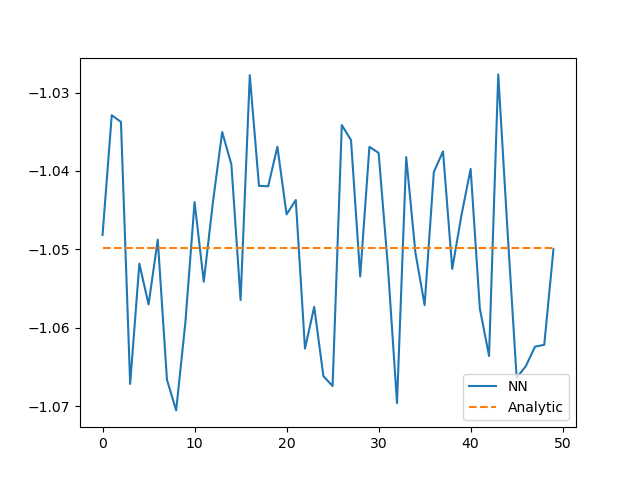}
    \caption*{Solution}
    \end{minipage}
        \begin{minipage}[c]{.38\linewidth}
    \includegraphics[width=0.90\linewidth]{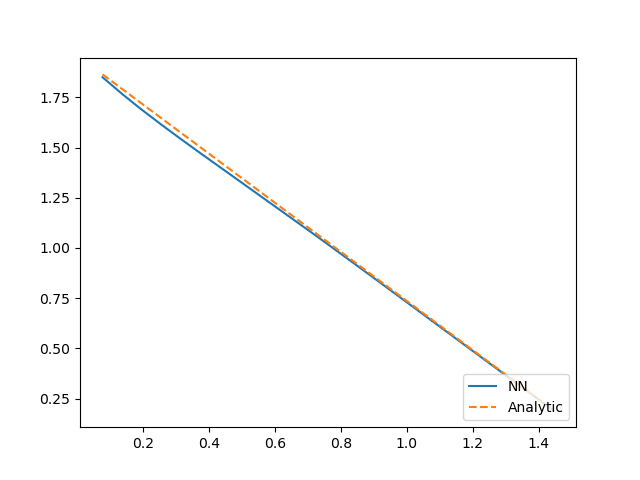}
    \caption*{Control}
    \end{minipage}
    \begin{minipage}[c]{.38\linewidth}
    \includegraphics[width=0.90\linewidth]{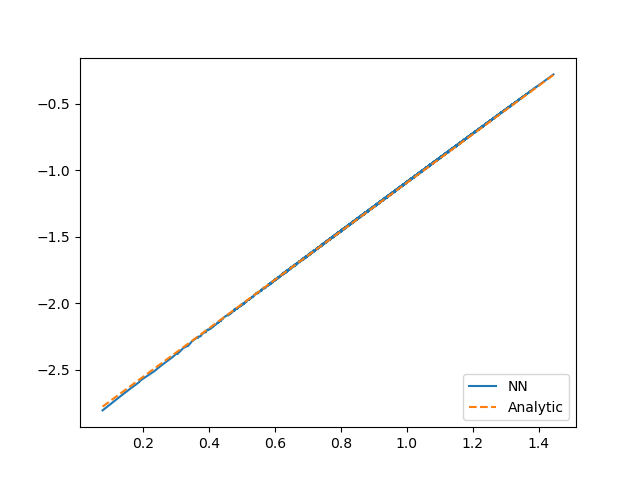}
    \caption*{Derivative}
    \end{minipage}
     \begin{minipage}[c]{.38\linewidth}
    \includegraphics[width=0.90\linewidth]{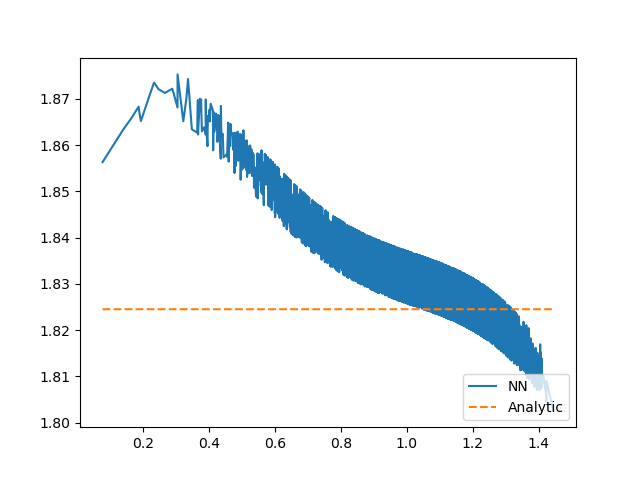}
    \caption*{Second order derivative}
    \end{minipage}
    \caption{\small{Solution and control obtained on the mean variance case at $\frac{T}{2}$ in dimension $300$ with $20$ time steps comparing analytic solution to the calculated one (NN). Truncation factor: $0.999$.}}
    \label{fig:solAnalMV300}
\end{figure}

Moreover, we test the accuracy of the control approximation. We solve the PDE from $T/2$ to $T$ starting  with  the optimal distribution of the wealth  at $T/2$, which is given by:
\begin{align} 
\log\Big(  \frac{X_{\frac{T}{2}}}{\eta(T/2)} - x_0 - e^{ \frac{RT}{2 \lambda}} \Big) &\sim \;  \mathcal{N}(0, \kappa(T/2)),  
\end{align}
with
\begin{align}
\eta(t)  &= \;  - e^{ \frac{R(T-t)}{2 \lambda}}, \quad  \kappa^2(t) \; = \;  \log\Big(\frac{e^{R(T-t)}(e^{RT}-e^{R(T-t)})}{4 \lambda^2 \eta(t)^2} +1 \Big),
\end{align}  
and we calculate the solution obtained  at date $T/2$ and the control obtained solving the PDE \eqref{PDE mean variance} that we can compare to the analytical solution given by  \eqref{eq:analMV} and \eqref{eq:optCont}.
After training, using $n_s=50$ samples of $X \in  (\R^N)^{n_s}$ following the law of $X_{\frac{T}{2}}$,   we calculate the control obtained for each sample in each of the dimension. After sorting $X$ in a one dimensional array,
We plot the result obtained  on Figures \ref{fig:solAnalMV100}--\ref{fig:solAnalMV300}. For the solution, the $x$-axis  corresponds to the sample number and the $y$-axis is the value of the estimated solution. For the other plots, the $x$-axis is the state space and the $y$-axis is the value of the corresponding function. 

\subsection{A min/max Linear quadratic mean-field control problem}

We consider a mean-field model in which the dynamics is linear and the running cost is quadratic in the position, the control and the expectation of the position. The terminal cost is encourages to be close to one of two targets. This type of model is inspired by the min-LQG problem of~\cite{salhab2015a}. More precisely, we consider the following controlled McKean-Vlasov dynamics 
\begin{align}
    \di X_t &= \; \big[  A X_t + \bar A \E[X_t] + B \alpha_t \big]  \ \di t +  \sigma \ \di W_t , 
\quad X_0 \sim \mu_0, 
\end{align}
where $\alpha$ $=$ $(\alpha_t)_t$  is the  control, and the agent aims to minimize  the functional cost
\begin{align} 
 \label{defV0-minLQ}
    J(\alpha) \; = \; \E \Big[ \int_0^T f(X_t,\E[X_t],\alpha_t) \ \di t + g(X_T) \Big] 
    &\;
    \rightarrow 
    \quad V_0 \; = \; \inf_{\alpha} J(\alpha), 
\end{align}
where the running and terminal costs are given by 
\begin{align}
    f(x,\bar x,a)  \; = \; \frac{1}{2} \left( Q x^2 + \bar Q (x - S \bar x)^2 + R a^2 \right), 
    & \quad g(x) \; = \; \min \left\{  |x - \xi_1 |^2, | x - \xi_2 |^2 \right\}, 
\end{align} 
for some non-negative constants $Q$, $\bar Q$, $S$, $R$, and two real numbers $\xi_1$ and $\xi_2$.  

The value function to the mean-field type control problem \eqref{defV0-minLQ} is solution to the Bellman (semi-linear PDE) equation \eqref{meanPDE} with $r$ $=$ $0$, and 
\begin{align}
    h(x,\mu,z,\gamma) 
    & = \; \inf_{a \in \R} \big\{  \big[  A x + \bar A  \E_\mu[\xi] + B a \big] z +  \frac{1}{2} \left( Q x^2 + \bar Q (x - S \E_\mu[\xi])^2 + R a^2 \right) \big\} + \frac{\sigma^2}{2}  \gamma   \\
    & = \;  \big[ A x + \bar A  \E_\mu[\xi] \big] z - \frac{B^2}{2R} z^2 +  \frac{1}{2} \left( Q x^2 + \bar Q (x - S \E_\mu[\xi])^2  \right) + \frac{\sigma^2}{2}  \gamma ,
\end{align}
where the minimizer in the above $\inf$ is given by $a = - \frac{B}{R} z$,  and the terminal condition $\Gc(\mu)$ $=$ $\E_\mu \min \left\{  |\xi - \xi_1 |^2, | \xi - \xi_2 |^2 \right\}$ is  the expected minimal distance to one of the  targets $\xi_1,\xi_2$.

For the sake of illustration, we present several test cases.

The targets are at $\xi_1 = 0.25$ and $\xi_2 = 1.75$. Here we used $A = \bar A = 0$, $B = 1$, $Q = 0, \bar Q = S = R = 1$, and a time horizon $T = 0.5$. The initial distribution $\mu_0$  is a Gaussian $\mathcal{N}(x_0, \vartheta_0^2)$. We consider the following test cases:
\begin{enumerate}
    \item $\sigma=0.3, x_0 = 1, \vartheta_0 = 0.2$,
    \item $\sigma=0.5, x_0 = 0.625, \vartheta_0 = \sqrt{0.2}$,
    \item $\sigma=0.3, x_0 = 0.625, \vartheta_0 = \sqrt{0.2}$,
    \item $\sigma=0.3, x_0 = 0.625, \vartheta_0 = \sqrt{0.4}$.
\end{enumerate}
References are given in table \ref{tab:MinLQCRef}: they are calculated by the PDE method in~\cite{MR2679575} (in the context of mean field games; see~\cite{MR3392611} for the adaptation to the PDE system arising in mean field control) with step size in space and time of size $10^{-3}$, and the neural network method referred to as Algorithm~1 in~\cite{carlau19} with $N=10000$ and $N_T$ $=$ $50$.
\begin{table}[H]
    \centering
    \begin{tabular}{|c|c|c|}
        \hline
  Case &  Benchmark & Global \\
        \hline
     1  &  $0.2256$   & $0.2273(0.004)$ \\ \hline 
     2 &  $0.2085$   & $0.2098(0.006)$ \\ \hline 
     3 &  $0.1734$   & $0.1742(0.005)$ \\ \hline 
     4 &  $0.2276$   & $0.2300(0.009)$ \\ \hline 
    \end{tabular}
    \caption{  \small{Min-LQC example reference solutions : benchmark solution estimated  by finite difference scheme and  Algorithm~1 in~\cite{carlau19} with  $N$ $=$ $10000$, $N_T$ $=$ $50$, 10 neurons and  3 hidden layers, tanh activation function, average on 10 runs.}}
    \label{tab:MinLQCRef}
\end{table}
In table \ref{tab:MinLQCNew}, we give the results obtained with different time discretization and dimension for the DPBD scheme  using ReLU activation functions with a  DeepSet network for $\Uc$ and a second AD-DeepSet  network  to estimate  $\Zc$. Results are very good except for test case 1 where a small bias appears.

\begin{small}
\begin{table}[H]
    \centering
    \begin{tabular}{|c|c|c|c|c|}
        \hline
  Case &   $N=100,N_T=30$& $N=100,N_T=60$& $N=500,N_T=30$ & $N=500,N_T=60$\\
        \hline
    1  &    $0.2370 (0.013)$ &  $0.2382 (0.012)$    & $0.2446   (0.013)$ & $0.2495 (0.09)$\\ \hline  
     2 &   $0.2088  (0.002)$ &  $0.2092 (0.001)$  & $0.2106  (0.003)$ & $0.2105 (0.003)$\\ \hline    
     3 &  $0.1774  (0.007)$ &  $0.1784 (0.005)$ & $0.1819 (0.005)$ & $0.1785 (0.008)$\\ \hline    
     4 &  $0.2279 (0.005)$ &  $0.2264 (0.005)$  & $0.2292 (0.006)$& $0.2274 (0.006)$\\ \hline  
    \end{tabular}
    \caption{  \small{Min-LQC example with DPBD scheme  using ReLU activation functions with a  DeepSet network for $\Uc$ and a second AD-DeepSet  network  to estimate  $\Zc$, average on 10 runs, standard deviation in parenthesis.}}
    \label{tab:MinLQCNew}
\end{table}
\end{small}

In table \ref{tab:MinLQCNew1Net}, we give the same results using a single network. Here the results are very good for all test cases.  Using two networks, the algorithm certainly face difficulties to approximate the derivatives near maturities which is not required using a single network.

\begin{table}[H]
    \centering
    \begin{tabular}{|c|c|c|c|c|}
        \hline
  Case &   $N=100,N_T=30$& $N=100,N_T=60$& $N=500,N_T=30$ & $N=500,N_T=60$\\
        \hline
    1  & $0.2289 (0.0006)$    & $0.2271 (0.001)$    & $0.2290 (0.0004)$ &  $0.2271  (0.0008)$\\ \hline  
     2 & $0.2083 (0.0008)$  & $0.2086 (0.0007)$ & $0.2097 (0.0008)$ & $0.2089  (0.0004)$ \\ \hline    
     3 & $0.1740 (0.001)$ & $0.1740 (0.001)$  & $0.1742 (0.0004)$ & $0.1729 (0.0007)$\\ \hline    
     4 & $0.2276 (0.001)$  & $0.2310 (0.003)$ & $0.2282 (0.0008)$&  $0.2278 (0.001)$\\ \hline  
    \end{tabular}
    \caption{  \small{Min-LQC example  with DBDP scheme  using ReLU activation functions and a single  DeepSet network for $\Uc$ which is differentiated to approximate  $\Zc$, average on 10 run, standard deviation in parenthesis.}}
    \label{tab:MinLQCNew1Net}
\end{table}

\begin{figure}[H]
    \centering
    \centering
    \begin{minipage}[c]{.24\linewidth}
    \includegraphics[width=\linewidth]{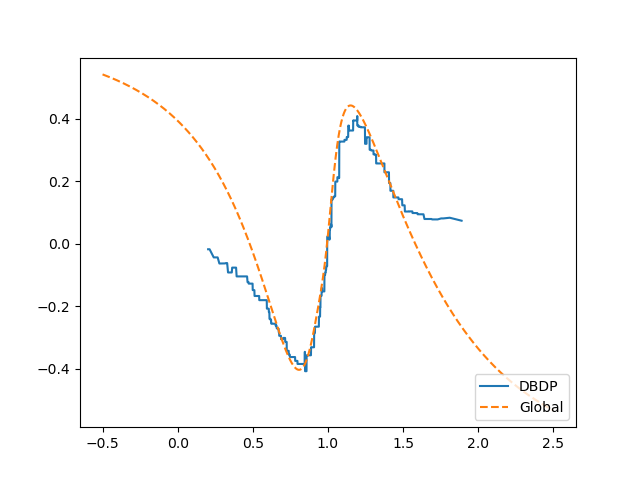}
    \caption*{Case 1}
    \end{minipage}
    \begin{minipage}[c]{.24\linewidth}
    \includegraphics[width=\linewidth]{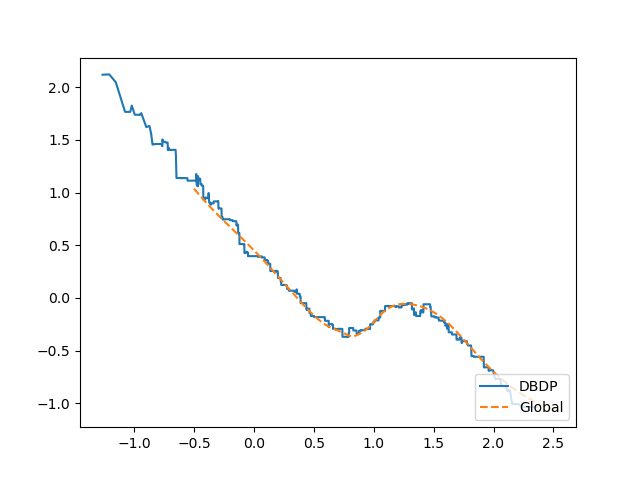}
    \caption*{Case 2}
    \end{minipage}
    \begin{minipage}[c]{.24\linewidth}
    \includegraphics[width=\linewidth]{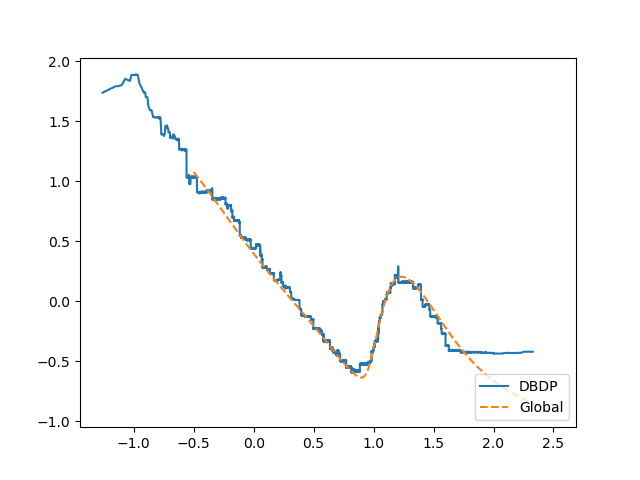}
    \caption*{Case 3}
    \end{minipage}
    \begin{minipage}[c]{.24\linewidth}
    \includegraphics[width=\linewidth]{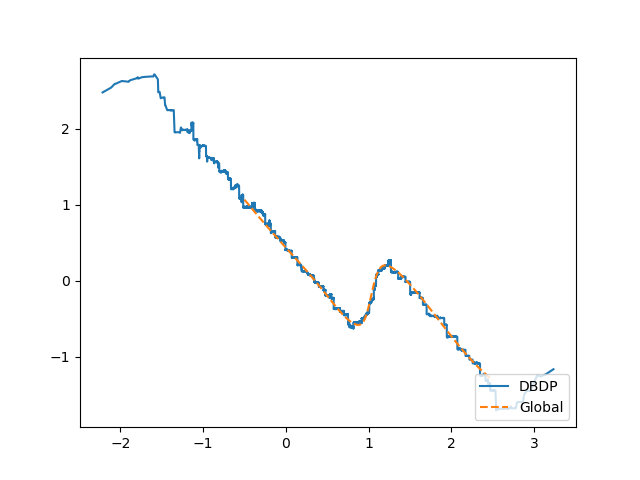}
    \caption*{Case 4}
    \end{minipage}
    \caption{Control calculated at $t=0$ for Min-LQC examples: comparison DBDP using a single DeepSet network with $N_T=50$, $N=500$ and global approximation. }
    \label{fig:solMinLQT0}
\end{figure}

\small

\printbibliography

\end{document}